\documentclass[twoside,11pt]{article}

\usepackage{amssymb, amsthm, amsmath}
\usepackage{graphicx}
\usepackage{hyperref}

\newtheorem{thm}{Theorem}[section]

\newtheorem{lem}[thm]{Lemma}
\newtheorem{prop}[thm]{Proposition}
\theoremstyle{definition}
\newtheorem{defn}[thm]{Definition}
\theoremstyle{remark}
\newtheorem{rem}[thm]{Remark}
\numberwithin{equation}{section}
\allowdisplaybreaks

\begin{document}
	\pagestyle{myheadings} \markboth{ \rm \centerline {Mateusz Kubiak and Bogdan Szal}} 
{\rm \centerline { }}

\begin{titlepage}
	\title{\bf {A sufficient condition for uniform convergence of double sine series with $p$-bounded variation coefficients}}
	
	\author {Mateusz Kubiak and Bogdan Szal \\
		{\small University of Zielona G\'{o}ra,}\\
		{\small Faculty of Mathematics, Computer Science and Econometrics,}\\
		{\small 65-516 Zielona G\'{o}ra, ul. Szafrana 4a, Poland} \\ 
		{\small e-mail: M.Kubiak@wmie.uz.zgora.pl, B.Szal @wmie.uz.zgora.pl}}
	
	
\end{titlepage}

\date{}
\maketitle
\begin{abstract}
	In the present paper we will introduce a new class of double sequences called $DGM\left(p, \alpha,\beta,\gamma,r\right) ,$ which is the generalization of a class considered by Szal and Duzinkiewicz. Moreover, we obtained in this note a sufficient condition for the uniform convergence of double sine series with coefficients belonging to this class.
\end{abstract}


	
\begin{section}{Introduction}
	\subsection{Uniform convergence of single sine series}
	It is well known that there is a great number of interesting results in
	Fourier analysis established by assuming monotonicity of Fourier coefficients. The
	following classical convergence result can be found in many monographs (see for example \cite{Chaundy} and \cite{Zygmund}).
	
	\begin{thm}\label{pierwsze twierdzenie}
		Suppose that $b_{n}\geq b_{n+1}$ and $b_{n}\rightarrow 0$. Then a necessary
		and sufficient condition for the uniform convergence of the series%
		\begin{equation}
			\sum\limits_{n=1}^{\infty }b_{n}\sin nx  \label{1}
		\end{equation}%
		is 
		\begin{equation}
			nb_{n}\rightarrow 0\text{ as n }\rightarrow \infty. \label{1.2}
		\end{equation}
	\end{thm}
	Several classes of sequences have been introduced to generalize Theorem \ref{pierwsze twierdzenie}. These classes are larger than the class of monotone sequences and contain sequences of complex numbers as well. The latest class is the following (see \cite{Kubiak}): 
	\begin{equation*}
		GM\left( p,\beta ,r\right) =\left\{ \left\{ a_{n}\right\} _{n=1}^{\infty
		}\subset \mathbb{C}:\exists _{C>0}\text{ }\forall _{m\in \mathbb{N}}\text{ }%
		\left\{ \sum\limits_{n=m}^{2m-1}\left\vert \Delta _{r}a_{n}\right\vert
		^{p}\right\} ^{\frac{1}{p}}\leq C\beta _{m}\right\} ,
	\end{equation*}%
	where $\Delta_{r} a_k = a_k - a_{k+r}$ for $ r \in \mathbb{N}$, $p>0$, $C$ is a positive constant depending only on $\left\{ a_{k}\right\} _{k=1}^{\infty }
	$ and a sequence $\left\{ \beta_k \right\}_{k=1}^{\infty} \subset \mathbb{R}_{+}$.  
	
	The classes: $GM(p,\beta,1)$ in \cite{Tikhonov1_new} and \cite{Liflyand_Tikhonov}, $GM(1,\beta,r)$ in  \cite{Szal} and \cite{Szal1}, and $GM(1,\beta,1)$ in \cite{Tikhonov1}, \cite{Tikhonov2}, \cite{Tikhonov3} and \cite{Leindler2} were defined, respectively.
	
	Theorem \ref{pierwsze twierdzenie} was generalized by supposing that the coefficients of (\ref{1}) belong to the class $GM(1,\dot{\beta},1)$, where $\dot{\beta} = \dot{\beta}_n = n^{-1} \sum\limits_{k=[n/\lambda]}^{[\lambda n]}\left\vert b_{k}\right\vert (\lambda > 1)$ and [$\cdot{}$] means integer part (see \cite{dyachenko}). K\'orus proved in \cite{Korus_new} that $MVBVS \subsetneq SBVS \subsetneq SBVS_2$, and a series (\ref{1}) with complex coefficients from the class $SBVS_2$ is uniformly convergent if (\ref{1.2}) is satisfied, where $MVBVS \equiv GM(1,\dot{\beta},1)$ with $\lambda \geq 2$ (see \cite{zhou}), $SBVS \equiv GM(1,\ddot{\beta},1)$ with $\ddot{\beta} = \ddot{\beta}_n = n^{-1} \underset{m \geq [n/\lambda]}{sup} \sum\limits_{k=m}^{2m}\left\vert c_k \right\vert (\lambda \geq 2)$ and $SBVS_2 \equiv GM(1,\dddot{\beta},1)$ with $\dddot{\beta} = \dddot{\beta}_n = n^{-1} \underset{m\geq b(n)}{sup} \sum\limits_{k=m}^{2m}\left\vert c_k \right\vert$ (where the sequence $\left\{b(k)\right\}_{k=1}^{\infty} \subset \mathbb{R}_+$ tends monotonically to infinity), respectively. In \cite{Szal} it is showed that the class $GM(1,\dot{\beta},2)$ is larger than $MVBVS$ and a series (\ref{1}) with coefficients from $GM(1,\dot{\beta},2)$ is uniformly convergent if (\ref{1.2}) holds.
	
	\subsection{Uniform convergence of double sine series}
	First, we introduce some definitions and notations.
	A double series
	\begin{equation*}
	\sum\limits_{j=1}^{\infty }\sum\limits_{k=1}^{\infty }a_{jk}
	\end{equation*}
	of complex numbers converges regularly if the sums
	\begin{equation*}
	\sum\limits_{j=1}^{m}\sum\limits_{k=1}^{n }a_{jk}
	\end{equation*}
	converge to finite number as $m$ and $n$ tend to infinity independently of each other, and moreover both the column and row series
	\begin{equation*}
	\sum\limits_{j=1}^{\infty}a_{jn},\text{ \ }n=1,2,...,\qquad
	\sum\limits_{k=1}^{\infty}a_{mk},\text{ \ }m=1,2,...,
	\end{equation*}%
	are convergent. Equivalently, for any $\epsilon >0$ there exists a positive number $m_0 = m_0(\epsilon)$ such that
	\begin{equation*}
	\left\vert \sum\limits_{j=1}^{M}\sum\limits_{k=1}^{N}a_{jk} \right\vert < \epsilon
	\end{equation*}
	for any $m,n,M,N$ with $m+n> m_0, 1\leq m \leq M$ and $1\leq n\leq N$.
	
	A monotonically decreasing double sequence of real nonnegative numbers $\left\{a_{jk} \right\}_{j,k =1}^{\infty}$ is a sequence such that
	\begin{equation*}
	\Delta_{10}a_{jk} \geq 0, \hspace{10mm}  \Delta_{01}a_{jk} \geq 0, \hspace{10mm} \Delta_{11}a_{jk} \geq 0, \hspace{10mm} j,k=1,2,...,
	\end{equation*}
	where for $r \in \mathbb{N}$
	\begin{equation*}
	\Delta_{r0}a_{jk} := a_{jk} - a_{j+r,k}, \hspace{15mm} \Delta_{0r}a_{jk} := a_{jk} - a_{j,k+r},
	\end{equation*}
	\begin{equation*}
	\Delta_{rr}a_{jk} := \Delta_{r0}\left( \Delta_{0r}a_{jk}\right) = a_{j,k} - a_{j+r,k} - a_{j,k+r} + a_{j+r,k+r}.
	\end{equation*}
	The class of all monotonically decreasing double sequences we will denote by $M$. Let $\left\{c_{jk} \right\}_{j,k =1}^{\infty}$ be a double sequence of complex numbers. Consider the double sine series
	\begin{equation}
	\sum\limits_{j=1}^{\infty }\sum\limits_{k=1}^{\infty }c_{jk} \sin jx \sin ky. \label{def_poswojny_szereg} 
	\end{equation}
	
	The following two-dimensional extension of Theorem \ref{pierwsze twierdzenie} is due to \v Zak and \v Sneider \cite{Zak}.
	\begin{thm}\label{twierdzenie_podstawowe_zak_sneider}
		If $\left\{c_{jk}\right\}_{j,k=1}^{\infty}\subset \mathbb{R}_{+}$ is a monotonically decreasing double sequence, then the series (\ref{def_poswojny_szereg}) is uniformly convergent in $(x,y)$ if and only if
		\begin{equation}\label{jkcjk_dozy_do_zera}
		jkc_{jk} \rightarrow 0\text{ as }j+k \rightarrow \infty.
		\end{equation}
	\end{thm}
	In 2009 K\'orus and M\'oricz defined a new class of sequence in the following way (see \cite{Korus_moricz}):
	\begin{defn}
		A double sequence $c :=\left\{ c_{jk}\right\}_{j,k=1}^{\infty} \subset \mathbb{C}$ belongs to the class $MVBVDS$ (called Mean Value Bounded Variation Double Sequences), if there exist a positive constant $C$ and an integer $\lambda \geq 2$, depending only on $\left\{ c_{jk}\right\}_{j,k=1}^{\infty}$, such that:
		\begin{equation}\label{eqn: 1_7}
		\sum\limits_{j=m}^{2m-1}\left\vert \Delta_{10}c_{jn}\right\vert \leq \frac{C}{m} \sum\limits_{j=[\lambda^{-1}m]}^{[\lambda m]} |c_{jn}| \hspace{10mm} m\geq \lambda, n = 1,2,...,
		\end{equation}
		\begin{equation}\label{eqn: 1_8}
		\sum\limits_{k=n}^{2n-1}\left\vert \Delta_{01}c_{mk}\right\vert \leq \frac{C}{n}  \sum\limits_{k=[\lambda^{-1}n]}^{[\lambda n] } |c_{mk}| \hspace{10mm} n\geq \lambda, m =  1,2,...,
		\end{equation}
		\begin{equation}\label{eqn: 1_9}
		\sum\limits_{j=m}^{2m-1}\sum\limits_{k=n}^{2n-1}\left\vert \Delta_{11}c_{jk}\right\vert \leq\frac{C}{nm} \sum\limits_{j=[\lambda^{-1}m]}^{[\lambda m]}\sum\limits_{k=\lambda^{-1}n}^{\lambda n } |c_{jk}| \hspace{10mm} m,n\geq \lambda.
		\end{equation}
	\end{defn}
	
	It is clear that $M \subsetneq MVBVDS$. In the paper \cite{Korus_moricz} K\'orus and M\'oricz proved the following theorem:
	\begin{thm}
		If $\left\{c_{jk} \right\}_{j,k =1}^{\infty} \subset \mathbb{C}$ belongs to $MVBVDS$ and  $(\ref{jkcjk_dozy_do_zera})$ holds, then the series (\ref{def_poswojny_szereg}) is uniformly regularly convergent in $(x,y)$.
	\end{thm}

	Next, K\'orus defined in \cite{Korus} new classes of sequences as follows:
	\begin{defn}(\cite{Korus})
		A double sequence $c :=\left\{ c_{jk}\right\}_{j,k=1}^{\infty} \subset \mathbb{C}$ belongs to the class $SBVDS_1$ if there exist constant $C$, an integer $\lambda \geq 2$ and sequences $\left\{ b_{1}(l)\right\}_{l=1}^{\infty}$, $\left\{ b_{2}(l)\right\}_{l=1}^{\infty}$, $\left\{ b_{3}(l)\right\}_{l=1}^{\infty}$, each one converges (not necessarily monotonically) to infinity, all of them depending only on $\left\{ c_{jk}\right\}_{j,k=1}^{\infty}$, such that:
		\begin{equation}\label{eqn: 1_10}
		\sum\limits_{j=m}^{2m-1}\left\vert \Delta_{10}c_{jn}\right\vert \leq \frac{C}{m} 
		\left( \underset{b_1(m)\leq M\leq \lambda b_1(m)}{\text{max}} \sum\limits_{j=M}^{2M} |c_{jn}|\right), m \geq \lambda, n \geq 1,
		\end{equation}
		\begin{equation}\label{eqn: 1_11}
		\sum\limits_{k=n}^{2n-1}\left\vert \Delta_{01}c_{mk}\right\vert \leq \frac{C}{n}  
		\left( \underset{b_2(n)\leq N\leq \lambda b_2(n)}{\text{max}} \sum\limits_{k=N}^{2N} |c_{mk}|\right), m \geq 1, n \geq \lambda,
		\end{equation}
		\begin{equation}\label{eqn: 1_12}
		\sum\limits_{j=m}^{2m-1}\sum\limits_{k=n}^{2n-1}\left\vert \Delta_{11}c_{jk}\right\vert \leq\frac{C}{nm} 
		\left( \underset{M+N\geq b_3(m+n)}{\text{sup}} \sum\limits_{j=M}^{2M}\sum\limits_{k=N}^{2N} |c_{jk}| \right),
		m,n\geq \lambda.
		\end{equation}
	\end{defn} 
	
	\begin{defn}
		(\cite{Korus}) A double sequence $\{c_{jk}\}_{j,k=1}^{\infty }\subset \mathbb{C}
		$ belongs to the class $SBVDS_{2}$, if there exist a positive constant $%
		\mathcal{C}$, an integer $\lambda \geq 1$ and a sequence $\{b(l)\}_{l=1}^{
			\infty }$ tending monotonically to infinity, depending only on $\{c_{jk}\}_{j,k=1}^\infty$,
		for which:
		\begin{equation}\label{supremum_1_0}
			\sum_{j=m}^{2m-1}|\Delta _{10}c_{jn}|\leq \frac{\mathcal{C}}{m}\Bigg(
			\sup_{M\geq b(m)}\sum_{j=M}^{2M}|{c_{jn}}|\Bigg) ,\ m\geq \lambda ,\ n\geq 1,
		\end{equation} 
		\begin{equation}\label{supremum_0_1}
			\sum_{k=n}^{2n-1}|\Delta _{01}c_{mk}|\leq  \frac{\mathcal{C}}{n}\Bigg(
			\sup_{N\geq b(n)}\sum_{k=N}^{2N}|{c_{mk}}|\Bigg) ,\ n\geq \lambda ,\ m\geq 1,
		\end{equation}
		\begin{equation}\label{supremum_1_1}
		\sum_{j=m}^{2m-1}\sum_{k=n}^{2n-1}|\Delta _{11}c_{jk}|\leq \frac{
			\mathcal{C}}{mn}\Bigg( \sup_{{M+N}\geq b(m+n)}\sum_{j=M}^{2M}\sum_{k=N}^{2N}{
			|c_{jk}|}\Bigg) ,\ m,n\geq \lambda \mathbf{.}
		\end{equation}
	\end{defn}

	In the same paper  K\'orus proved the following theorem:
	\begin{thm}\label{twierdzenie_svbvds_2}(\cite{Korus})
		If $\{c_{jk}\}_{j,k=1}^{\infty }\subset \mathbb{C}$ belongs to the class $SBVDS_{2}$ and \eqref{jkcjk_dozy_do_zera} holds, then the series \eqref{def_poswojny_szereg} is uniformly regularly convergent in $(x,y)$.
		\label{T.4}
	\end{thm}
	He showed also that $MVBVDS \subsetneq SBVDS_1 \subsetneq SBVDS_2 $.
	
	Next, Szal and Duzinkiewicz defined a new class of sequence in the following way (see \cite{Szal_Duzinkiewicz}):
	\begin{defn}
		A double sequence $c :=\left\{ c_{jk}\right\}_{j,k=1}^{\infty} \subset \mathbb{C}$ belongs to class $DGM(\alpha,\beta, \gamma,r)$ (called Double General monotone), if there exist a positive constant $C$ and an integer $\lambda \geq1$, depending only on $\left( c_{jk}\right)_{j,k=1}^{\infty}$, for which:
		\begin{equation*}
		\sum\limits_{j=m}^{2m-1}\left\vert \Delta_{r0}c_{jn}\right\vert \leq C \alpha_{mn} \hspace{10mm} m\geq \lambda, n\geq 1,
		\end{equation*}
		\begin{equation*}
		\sum\limits_{k=n}^{2n-1}\left\vert \Delta_{0r}c_{mk}\right\vert \leq C \beta_{mn} \hspace{10mm} n\geq \lambda, m\geq 1,
		\end{equation*}
		\begin{equation*}
		\sum\limits_{j=m}^{2m-1}\sum\limits_{k=n}^{2n-1}\left\vert \Delta_{rr}c_{jk}\right\vert \leq C \gamma_{mn} \hspace{10mm} m,n\geq \lambda
		\end{equation*}
		hold, where $\alpha :=\left\{ \alpha_{jk}\right\}_{j,k=1}^{\infty}$, $\beta :=\left\{ \beta_{jk}\right\}_{j,k=1}^{\infty}$, $\gamma :=\left\{ \gamma_{jk}\right\}_{j,k=1}^{\infty}$ are nonnegative double sequences and $r\in \mathbb{N}$.
	\end{defn}
	
	Using this definition with $r=1$, we have:
	\begin{itemize}
		\item $MVBVDS$ $\equiv$ $DGM({}_1\alpha,{}_1 \beta,{}_1\gamma,1)$, where $\left\{ {}_1\alpha \right\}$, $\left\{{}_1\beta \right\}$ and $\left\{{}_1\gamma \right\}$ are the sequences defined by the formulas on the right sides of the inequalities \eqref{eqn: 1_7}, \eqref{eqn: 1_8} and \eqref{eqn: 1_9}, respectively;
		
		\item $SVBVDS_1$ $\equiv$ $DGM({}_2\alpha,{}_2\beta,{}_2\gamma,1)$, where $\left\{ {}_2\alpha \right\}$, $\left\{{}_2\beta \right\}$ and $\left\{{}_2\gamma \right\}$ are the sequences defined by the formulas on the right sides of the inequalities \eqref{eqn: 1_10}, \eqref{eqn: 1_11} and \eqref{eqn: 1_12}, respectively.
		
		\item $SVBVDS_2$ $\equiv$ $DGM({}_3\alpha,{}_3\beta,{}_3\gamma,1)$, where $\left\{ {}_3\alpha \right\}$, $\left\{{}_3\beta \right\}$ and $\left\{{}_3\gamma \right\}$ are the sequences defined by the formulas on the right sides of the inequalities \eqref{supremum_1_0}, \eqref{supremum_0_1} and \eqref{supremum_1_1}, respectively.
	\end{itemize} 
	
	Szal and Duzinkiewicz proved in \cite{Szal_Duzinkiewicz} the following theorem:
	\begin{thm}\label{twierdzenie_Szal_Duzinkiewicz}
		If a double sequence $\left\{c_{jk} \right\}_{j,k =1}^{\infty} \subset \mathbb{C}$ belongs to $DGM({}_2\alpha,{}_2\beta,{}_2\gamma,2)$ and (\ref{jkcjk_dozy_do_zera}) holds, then the series \eqref{def_poswojny_szereg} is uniformly regularly convergent in $(x,y)$.
	\end{thm}
	They showed also in this paper that $DGM({}_2\alpha,{}_2\beta,{}_2\gamma,1) \subsetneq DGB({}_2\alpha,{}_2\beta,{}_2\gamma,2)$.

	The latest class of sequences was defined by Kubiak and Szal in \cite{Kubiak_Szal_podwojny} as follows:
	
	\begin{defn}
		A double sequence $c :=\left\{ c_{jk}\right\}_{j,k=1}^{\infty} \subset \mathbb{C}$ belongs to the class $DGM(p,\alpha,\beta, \gamma,r)$ if there exist a positive constant $C$ and an integer $\lambda \geq 1$, depending only on $\left\{ c_{jk}\right\}_{j,k=1}^{\infty}$, for which:
		\begin{equation*}
		\left(\sum\limits_{j=m}^{2m-1}\left\vert \Delta_{r0}c_{jn}\right\vert^p\right)^{\frac{1}{p}} \leq C \alpha_m \hspace{10mm} m\geq \lambda, n\geq 1,
		\end{equation*}
		\begin{equation*}
		\left(\sum\limits_{k=n}^{2n-1}\left\vert \Delta_{0r}c_{mk}\right\vert^p\right)^{\frac{1}{p}} \leq C \beta_n \hspace{10mm} n\geq \lambda, m\geq 1,
		\end{equation*}
		\begin{equation*}
		\left(\sum\limits_{j=m}^{2m-1}\sum\limits_{k=n}^{2n-1}\left\vert \Delta_{rr}c_{jk}\right\vert^p\right)^{\frac{1}{p}} \leq C \gamma_{mn} \hspace{10mm} m,n\geq \lambda
		\end{equation*}
		hold, where $\alpha :=\left\{ \alpha_{jk}\right\}_{j,k=1}^{\infty}$, $\beta :=\left\{ \beta_{jk}\right\}_{j,k=1}^{\infty}$, $\gamma :=\left\{ \gamma_{jk}\right\}_{j,k=1}^{\infty}$ are nonnegative double sequences, $r\in \mathbb{N}$ and $p >0$.
	\end{defn}
	It is clear that $DGM(1,\alpha,\beta,\gamma,r) = DGM(\alpha,\beta,\gamma,r)$. 
	
	Moreover Kubiak and Szal showed in the paper \cite{Kubiak_Szal_podwojny} the following properties:
	\begin{flalign*}
		& DGM(p_1,{}_2\alpha,{}_2\beta,{}_2\gamma,r) \subseteq DGM(p_2,{}_2\alpha,{}_2\beta,{}_2\gamma,r)\text{, where }0<p_1\leq p_2, r \in \mathbb{N},\\
		&
		DGM(p,{}_2\alpha,{}_2\beta,{}_2\gamma,r_1) \subseteq DGM(p,{}_2\alpha,{}_2\beta,{}_2\gamma,r_2)\text{, where }p \geq 1, r_1,r_2 \in \mathbb{N}\text{ and }r_1 | r_2.
	\end{flalign*}
	
	In this paper we generalize Theorem \ref{twierdzenie_podstawowe_zak_sneider} to the class $DGM(p,{}_2\alpha,{}_2\beta,{}_2\gamma,r)$ for $p>1$ and $r \in \mathbb{N}$ with some additional assumptions.
\end{section}

\begin{section}{Main Results}
	
	\begin{thm}\label{twierdzenie_maximum}
		Let a sequence $\left\{ c_{jk} \right\}_{j,k=0}^{\infty} \subset \mathbb{C}$ belong to $DGM\left(p,{}_2\alpha, {}_2\beta,{}_2\gamma,r\right) $,
		where $p > 1$ and $r \in \mathbb{N}$. Suppose that the series
		\begin{equation}\label{zbieznosc_w_punktach}
			\sum\limits_{j=1}^{\infty}\sum\limits_{k=1}^{\infty} c_{jk} \sin \left(\frac{2l_1 \pi }{r}j\right)  \sin \left(\frac{2l_2 \pi }{r} k\right)\text{, where }r \geq 3
		\end{equation}
		are regularly convergent for all $l_1,l_2 \in \{ 1,...,[\frac{r}{2}]-1\}$ when $r$ is an even number and $l_1,l_2 \in \{  1,...,[\frac{r}{2}]\}$ when $r$ is an odd number, and
		\begin{equation}\label{dodatkowy_warunek_dla_max}
			\text{ }\underset{j\geq m}{sup} \text{ } j \ln j\sum\limits_{k=n}^{\infty}  \left\vert c_{jk} \right\vert \rightarrow 0,
		\end{equation}
		\begin{equation}\label{dodatkowy_warunek_dla_max_k}
			\text{ }\underset{k\geq n}{sup} \text{ } k \ln k\sum\limits_{j=m}^{\infty}  \left\vert c_{jk} \right\vert \rightarrow 0
		\end{equation}
		as $m+n \rightarrow \infty$. If
		\begin{equation}\label{prof 13}
			mn\ln m\ln n \left\vert c_{mn}\right\vert \rightarrow 0\text{ as } %
			m+n\rightarrow \infty,
		\end{equation}
		then the series \eqref{def_poswojny_szereg} is uniformly regularly convergent in $(x,y)$.
	\end{thm}

	\begin{rem}
		We can notice that the condition \eqref{zbieznosc_w_punktach} is obviously satisfied for  $r = 1$ and $r = 2$.
	\end{rem}

	The following proposition shows that the condition \eqref{jkcjk_dozy_do_zera} does not guarantee the uniform convergence of the series \eqref{def_poswojny_szereg} when $\{c_{jk}\}_{j,k=0}^{\infty} \in DGM\left(p,{}_2\alpha, {}_2\beta,{}_2\gamma,3 \right)$ for $p>1$.
	\begin{prop}\label{proposition_rozbiezne_w_punkcie}
		There exist an $\left(x_0,y_0\right) \in \mathbb{R}^2$ and a sequence $\left\{ c_{jk} \right\}_{j,k=0}^\infty$ belonging to $DGM(p,{}_2\alpha, {}_2\beta, {}_2\gamma, 3)$ for $p>1$ with the properties $jkc_{jk} \rightarrow 0 $ as $n \rightarrow \infty$ such that $\left\{ c_{jk} \right\}_{j,k=0}^\infty$ does not belong to $DGM(1,{}_2\alpha, {}_2\beta, {}_2\gamma, 3)$, for which the series \eqref{def_poswojny_szereg} is divergent in $\left(x_0,y_0\right) $.
	\end{prop}
\end{section}
	
\begin{section}{Lemmas}
	
	Denote, for $r \in \mathbb{N}$ and $k=0,1,2...$ by
	\begin{equation*}
		\tilde{D}_{k,r}{(x)}=\frac{\cos \left(k+\frac{r}{2}\right)x}{2\sin \frac{rx}{2}}, \hspace{15mm} 
		D_{k,r}{(x)}=\frac{\sin \left(k+\frac{r}{2}\right)x}{2\sin \frac{rx}{2}}
	\end{equation*}
	the Dirichlet type kernels.

	\begin{lem} \text{(\cite{Szal}, \cite{Szal1}} \label{roznicowanie})
		Let $r,m ,n \in \mathbb{N}, l\in \mathbb{Z}$ and $\left\{a_k\right\}_{k=1}^{\infty} \subset \mathbb{C}$. If $x \neq \frac{2l\pi}{r}$, then for all $m \geq n$
		\begin{equation*}
		\sum\limits_{k=n}^{m} a_{k} sin(kx) = 
		- \sum\limits_{k=n}^{m} \Delta_{r} a_{k} \tilde{D}_{k,r}(x)
		+
		\sum\limits_{k=m+1}^{m+r} a_{k} \tilde{D}_{k,-r}(x) 
		-
		\sum\limits_{k=n}^{n+r-1} a_{k} \tilde{D}_{k,-r}(x),
		\end{equation*}
		where $\Delta_ra_k = a_k - a_{k+r}$.
	\end{lem}
	
	\begin{lem} \label{lemat_do_glownego_tw}
		Let $\left\{ a_{k} \right\}_{k=0}^{\infty} \subset \mathbb{C}$, $l\in \mathbb{Z}$ and $n,N \in \mathbb{N}$ be such that $n \leq N$.
		If $x \in (\frac{2l\pi}{r},\frac{(2l+1)\pi}{r})$, then
		\begin{equation}\label{eq_roznicowanie_pierwsza_polowka}
		\left\vert  \sum\limits_{k=n}^{N} a_{k} sin(kx) \right\vert \leq 
		\frac{\pi}{2\left(rx-2\pi l\right)} \left( \sum\limits_{k=n}^{N} \left\vert\Delta_{r}a_{k}\right\vert +
		\sum\limits_{k=N+1}^{N+r} \left\vert a_{k}\right\vert +
		\sum\limits_{k=n}^{n+r-1}\left\vert a_{k}\right\vert
		\right),
		\end{equation}
		and if $x \in (\frac{(2l+1)\pi}{r},\frac{(2l+2)\pi}{r})$, then
		\begin{equation}\label{eq_roznicowanie_druga_polowka}
		\left\vert  \sum\limits_{k=n}^{N} a_{k} sin(kx) \right\vert \leq 
		\frac{\pi}{2(l+1)\pi - rx} \left( \sum\limits_{k=n}^{N} \left\vert\Delta_{r}a_{k}\right\vert +
		\sum\limits_{j=N+1}^{N+r} \left\vert a_{k}\right\vert +
		\sum\limits_{j=n}^{n+r-1}\left\vert a_{k}\right\vert
		\right).
		\end{equation}

		\begin{proof}
			From Lemma \ref{roznicowanie} we have
			\begin{equation*}
			\left\vert  \sum\limits_{k=n}^{N} a_{k} sin(kx) \right\vert \leq 
			\sum\limits_{k=n}^{N} \left\vert \Delta_{r} a_{k} \right\vert \left\vert \tilde{D}_{k,r}(x)\right\vert
			+
			\sum\limits_{k=N+1}^{N+r} \left\vert a_{k} \right\vert \left\vert \tilde{D}_{k,-r}(x)\right\vert 
			+
			\sum\limits_{k=n}^{n+r-1} \left\vert a_{k} \right\vert \left\vert \tilde{D}_{k,-r}(x)\right\vert .
			\end{equation*}
			If $x \in (\frac{2l\pi}{r},\frac{(2l+1)\pi}{r})$, then using inequality $\left\vert \sin(\frac{rx}{2}) \right\vert  \geq \frac{rx}{\pi} -2l$ we obtain
			\begin{equation}\label{oszacowanie_jadro_diri}
			\left\vert \tilde{D}_{k,\pm r}(x) \right\vert= \left\vert \frac{cos(k \pm \frac{r}{2})x}{2sin(\pm \frac{rx}{2})} \right\vert
			\leq \frac{1}{ \left\vert 2sin(\pm \frac{rx}{2})\right\vert} 
			\leq  \frac{1}{2\left(\frac{rx}{\pi}-2l\right)}.
			\end{equation}
			From this we get
			\begin{equation*}
			\left\vert  \sum\limits_{k=n}^{N} a_{k} sin(kx) \right\vert 
			\leq 
			\frac{\pi}{2\left(rx-2l\pi\right)} \left(
			\sum\limits_{k=n}^{N} \left\vert \Delta_{r} a_{k} \right\vert
			+
			\sum\limits_{k=N+1}^{N+r} \left\vert a_{k} \right\vert 
			+
			\sum\limits_{k=n}^{n+r-1} \left\vert a_{k} \right\vert \right).
			\end{equation*}
			If $x \in (\frac{(2l+1)\pi}{r},\frac{(2l+2)\pi}{r})$, then using inequality $\left\vert \sin(\frac{rx}{2}) \right\vert \geq 2(l+1) - \frac{rx}{\pi} $ we obtain
			\begin{eqnarray*}
				\left\vert \tilde{D}_{k, \pm r}(x) \right\vert \leq \frac{1}{2\left( 2\left( l+1 \right) - \frac{rx}{\pi} \right)}
			\end{eqnarray*}
		    and
			\begin{equation*}
			\left\vert  \sum\limits_{k=n}^{N} a_{k} sin(kx) \right\vert 
			\leq 
			\frac{\pi}{2\left( 2(l+1) - \frac{rx}{\pi}\right)} \left(
			\sum\limits_{k=n}^{N} \left\vert \Delta_{r} a_{k} \right\vert
			+
			\sum\limits_{k=N+1}^{N+r} \left\vert a_{k} \right\vert 
			+
			\sum\limits_{k=n}^{n+r-1} \left\vert a_{k} \right\vert \right).
			\end{equation*}
		This ends the proof.
		\end{proof}
	\end{lem}

	\begin{lem}\label{lemat_do_warunku_z_logarytmem} (\cite{Kubiak_pojedynczy_dostateczny})
		Let $n, N \in \mathbb{N}$. Then for $p \geq 1$
		\begin{equation*}
		\int_{n+N^{\frac{1}{p}}}^{n+N} \frac{1}{k\ln k} dk \leq \ln p.
		\end{equation*}
		\begin{proof}
			This inequality is true for $p =1 $. Consider the function $f(p) = \left(n+N^{\frac{1}{p}}\right)^p$ for $p > 0$. We get:
			\begin{equation*}
				f'(p) = p\left(n+N^{\frac{1}{p}}\right)^{p-1} \frac{1}{p}N^{\frac{1}{p}-1} =  \left(n+N^{\frac{1}{p}}\right)^{p-1} N^{\frac{1}{p}-1} \geq 0 \text{ for all }p>0.
			\end{equation*}
			It means that the function is nondecreasing with respect to $p$. 
			\newline Thus:
			\begin{equation}
				n+N  = f(1) \leq f(p) = \left(n+N^{\frac{1}{p}}\right)^p \text{ for } p\geq 1.
			\end{equation}
			Hence we get that:
			\begin{equation}\label{nierowosclemat_z_calka}
				\ln\left(n+N\right) \leq \ln\left(n+N^{\frac{1}{p}}\right)^p.
			\end{equation}
			Therefore, integrating by substitution with $\ln k = t$ and using (\ref{nierowosclemat_z_calka}), we get
			\begin{equation*}
				\int_{n+N^{\frac{1}{p}}}^{n+N} \frac{1}{k \ln k } dk =
				\int_{\ln(n+N)^{\frac{1}{p}}}^{\ln(n+N)} \frac{1}{t} dt
				= \ln(\ln(n+N)) - \ln(\ln(n+N^\frac{1}{p}))
			\end{equation*}
			\begin{equation*}
				= \ln\left(\frac{\ln(n+N)}{\ln\left(n+N^{\frac{1}{p}}\right)}\right)
				= \ln\left(p \frac{\ln\left(n+N\right)}{\ln\left(n+N^{\frac{1}{p}}\right)^p}\right)
				\leq \ln p
			\end{equation*}
			and the proof is completed.
		\end{proof}
	\end{lem}

	\begin{lem}\label{zbieznosc_podwojny_bez_supremum}
		If $\left\{ c_{jk} \right\}_{j,k=0}^{\infty} \subset \mathbb{C}$, $\left( b\left(n\right) \right)_{n=1}^{\infty}$ is a nonnegative sequence such that $b(n) \nearrow \infty$ as $n \rightarrow \infty$, $p \geq 1$, \eqref{jkcjk_dozy_do_zera} and 
		\begin{equation*}
		\left(\sum\limits_{j=m}^{2m-1}\sum\limits_{k=m}^{2m-1}\left\vert \Delta_{rr}c_{jk}\right\vert^p\right)^{\frac{1}{p}} \leq  
		\frac{C}{nm} \underset{M + N \geq b(m+n)}{sup} \sum\limits_{k=M}^{2M}\sum\limits_{j=N}^{2N} |c_{jk}|, \text{ for $m,n \in \mathbb{N}$}
		\end{equation*}
		hold, then 
		\begin{equation}\label{lemat_mn_do_twierdzenia}
		m^{\frac{1}{p}}n^{\frac{1}{p}}\sum\limits_{j=m}^{\infty} \sum\limits_{k=n}^{\infty} \left\vert \Delta_{rr}c_{jk}\right\vert \rightarrow 0 \text{, as } m+n \rightarrow \infty.
		\end{equation}
		\begin{proof}
			Let $\varepsilon > 0$. Using the H\"{o}lder inequality with $p>1$ and \eqref{jkcjk_dozy_do_zera} we have
			\begin{flalign*}
			\sum\limits_{n=m}^{\infty}\sum\limits_{n=m}^{\infty} & \left\vert \Delta_{rr}c_{jk}\right\vert
			= \sum\limits_{r=0}^{\infty} \sum\limits_{s=0}^{\infty} \sum\limits_{j=2^rm}^{2^{r+1}m-1} \sum\limits_{k=2^sn}^{2^{s+1}n-1} \left\vert \Delta_{rr}c_{jk}\right\vert & \\
			& 
			\leq
			\sum\limits_{r=0}^{\infty} \sum\limits_{s=0}^{\infty} \left\{ \left(\sum\limits_{j=2^rm}^{2^{r+1}m-1} \sum\limits_{k=2^sn}^{2^{s+1}n-1} \left\vert \Delta_{rr}c_{jk}\right\vert^p \right)^{\frac{1}{p}} 
			\left(\sum\limits_{j=2^rm}^{2^{r+1}m-1} \sum\limits_{k=2^sn}^{2^{s+1}n-1}  1 \right)^{1-\frac{1}{p}} 
			\right\} \\
			&
			\leq
			\sum\limits_{r=0}^{\infty} \sum\limits_{s=0}^{\infty} \left\{ 
			\frac{C}{2^rm2^sn} \underset{M+N\geq b(2^rm+2^sn)}{sup} \sum\limits_{k=M}^{2M}\sum\limits_{j=N}^{2N} |c_{jk}|
			\left(2^{s}m 2^{r}n   \right)^{1-\frac{1}{p}} 
			\right\} \\
			&
			\leq
			\sum\limits_{r=0}^{\infty} \sum\limits_{s=0}^{\infty} \left\{ 
			\frac{C}{\left(2^rm2^sn\right)^{\frac{1}{p}} } \underset{M+N\geq b(2^rm+2^sn)}{sup} \sum\limits_{k=M}^{2M}\sum\limits_{j=N}^{2N} \frac{kj|c_{jk}|}{kj}
			\right\} \\
			&
			<
			\frac{C\varepsilon}{\left(mn\right)^{\frac{1}{p}}}\sum\limits_{r=0}^{\infty} \sum\limits_{s=0}^{\infty} \left\{ 
			\frac{1}{\left(2^r2^s\right)^{\frac{1}{p}}} \underset{M+N\geq b(2^rm+2^sn)}{sup}
			\sum\limits_{k=M}^{2M}\sum\limits_{j=N}^{2N} \frac{1}{kj}
			\right\} \\
			&
			\leq \frac{4C\varepsilon}{\left(mn\right)^{\frac{1}{p}}}\sum\limits_{r=0}^{\infty} \sum\limits_{s=0}^{\infty} \left\{ 
			\frac{1}{\left(2^r2^s\right)^{\frac{1}{p}}} \right\}
			\ll \frac{\varepsilon}{m^{\frac{1}{p}}n^{\frac{1}{p}}}.
			\end{flalign*}
			This implies that (\ref{lemat_mn_do_twierdzenia}) holds.
		\end{proof}
	\end{lem}

	\begin{lem}\label{zbieznosc_pojedynczych_z_supremum}
		Let $\left\{ c_{jk} \right\}_{j,k=0}^{\infty} \subset \mathbb{C}$, $\left( b\left(n\right) \right)_{n=1}^{\infty}$ be a nonnegative sequence such that $b(n) \nearrow \infty$ as $n \rightarrow \infty$, $\lambda \geq 2$ and  \eqref{jkcjk_dozy_do_zera} hold. If
		\begin{equation*}
		\left(\sum\limits_{j=m}^{2m-1}\left\vert \Delta_{r0}c_{jk}\right\vert^p\right)^{\frac{1}{p}} \leq  
		\frac{C}{m} \underset{b(m)\leq M\leq \lambda b(m)}{\text{max}} \sum\limits_{j=M}^{2M} |c_{jn}| \text{     for $m, k \in \mathbb{N}$} \text{ and } m \geq \lambda,  
		\end{equation*}
		then 
		\begin{equation}\label{lemat_m_do_twierdzenia_supremum}
		m^{\frac{1}{p}} \text{ }\underset{k\geq n}{sup} \text{ }k  \sum\limits_{j=n}^{\infty} \left\vert \Delta_{r0}c_{jk}\right\vert  \rightarrow 0\text{, as } m+n \rightarrow \infty.
		\end{equation}
		And if
		\begin{equation*}
		\left(\sum\limits_{k=n}^{2n-1}\left\vert \Delta_{0r}c_{jk}\right\vert^p\right)^{\frac{1}{p}} \leq  
		\frac{C}{n} \underset{b(n)\leq N\leq \lambda b(n)}{\text{max}}  \sum\limits_{k=N}^{2N} |c_{jk}| \text{     for $n, j \in \mathbb{N}$} \text{ and } n \geq \lambda,
		\end{equation*}
		then 
		\begin{equation}\label{lemat_n_do_twierdzenia_supremum}
		n^{\frac{1}{p}} \text{ }\underset{j\geq m}{sup} \text{ }j  \sum\limits_{k=n}^{\infty} \left\vert \Delta_{0r}c_{jk}\right\vert  \rightarrow 0\text{, as } m+n \rightarrow \infty.
		\end{equation}
		\begin{proof}
			Let $\varepsilon > 0$. Using the H\"{o}lder inequality with $p>1$ and \eqref{jkcjk_dozy_do_zera} we obtain
			\begin{flalign*}
			\underset{k\geq n}{sup} \text{ }k  \sum\limits_{k=n}^{\infty} & \left\vert \Delta_{r0}c_{jk}\right\vert
			=
			\underset{k\geq n}{sup} \text{ }k  \sum\limits_{r=0}^{\infty} \sum\limits_{j=2^rm}^{2^{r+1}m-1} \left\vert \Delta_{r0}c_{jk}\right\vert & \\
			&
			\leq
			\underset{k\geq n}{sup} \text{ } k \sum\limits_{r=0}^{\infty} \left(\left(\sum\limits_{j=2^rm}^{2^{r+1}m-1} \left\vert \Delta_{r0}c_{jk}\right\vert^p
			\right)^\frac{1}{p}
			\left(\sum\limits_{j=2^rm}^{2^{r+1}m-1} 1\right)^{1-\frac{1}{p}}\right) \\
			&
			\leq
			\underset{k\geq n}{sup} \text{ }k  \sum\limits_{r=0}^{\infty} 
			\left(
			\frac{C}{2^rm} \underset{b(m)\leq M\leq \lambda b(m)}{\text{max}}  \sum\limits_{j=M}^{2M} |c_{jk}|
			\left(2^{r}m\right)^{1-\frac{1}{p}}
			\right) \\
			&
			<
			\frac{C\varepsilon}{m^{\frac{1}{p}}}
			\underset{k\geq n}{sup} \text{ }k  \sum\limits_{r=0}^{\infty} 
			\left(
			\frac{1}{\left(2^r\right)^{\frac{1}{p}}} \underset{b(m)\leq M\leq \lambda b(m)}{\text{max}} \sum\limits_{j=M}^{2M} \frac{1}{kj}
			\right) 
			\leq
			\frac{2C\lambda\varepsilon}{m^{\frac{1}{p}}}  \sum\limits_{r=0}^{\infty} 
			\frac{1}{\left( 2^r\right)^{\frac{1}{p}} }
			\ll \frac{\varepsilon}{m^{\frac{1}{p}}}.
			\end{flalign*}
			This implies that (\ref{lemat_m_do_twierdzenia_supremum}) holds.
			Similarly as above, we can prove (\ref{lemat_n_do_twierdzenia_supremum}).
		\end{proof}
	\end{lem}

\end{section}

\begin{section}{Proofs}
	
\subsection{Proof of the theorem \ref{twierdzenie_maximum}}

	Analogously as in \cite[Theorem 12]{Kubiak_pojedynczy_dostateczny} we can show that the single series: 
	\begin{align}
	\sum_{j=1}^{\infty }c_{jn}\sin jx,\qquad n=1,2,\ldots ,\qquad
	\sum_{k=1}^{\infty }c_{mk}\sin ky,\qquad m=1,2,\ldots  \label{4.1}
	\end{align}
	are uniformly convergent since $\{c_{jn}\}_{j=1}^{\infty }\in GM(p,_{3}\beta,r)$ for any $n\in \mathbb{N}$ and $\{c_{mk}\}_{k=1}^{\infty }\in
	GM(p,_{3}\beta,r)$ for any $m\in \mathbb{N}$. Let $\epsilon >0$ be
	given, $p>1$ and $r \in \mathbb{N}$. We shall prove that for any $M \geq m > \eta (\varepsilon)$, $N \geq n > \eta (\varepsilon)$ and any $(x,y) \in \mathbb{R}^2$:
	\begin{equation}\label{zbieznosc_reszty_szeregu_podwojnego}
	\left\vert \sum\limits_{j=m}^{M}\sum\limits_{k=n}^{N} c_{jk} \sin(jx)\sin(ky) \right\vert \ll \varepsilon,
	\end{equation}
	where $\eta = \eta(\varepsilon) > \lambda$ is a natural number such that for any $m, n > \eta$
	\begin{equation}\label{lematy_zbieznosc_w_twierdzeniu}
	\begin{split}
		& \left\vert \sum\limits_{j=m}^{M}\sum\limits_{k=n}^{N} c_{jk} \sin(j\frac{2l_1\pi}{r})\sin(k\frac{2l_2\pi}{r})\right\vert < \varepsilon ,\\
		& m^{\frac{1}{p}} \text{ }\underset{k\geq n}{sup} \text{ }k  \sum\limits_{j=n}^{\infty} \left\vert \Delta_{r0}c_{jk}\right\vert < \varepsilon, \qquad n^{\frac{1}{p}} \text{ }\underset{j\geq m}{sup} \text{ }j  \sum\limits_{k=n}^{\infty} \left\vert \Delta_{r0}c_{jk}\right\vert < \varepsilon, \\
		& m^{\frac{1}{p}}n^{\frac{1}{p}}\sum\limits_{j=m}^{\infty} \sum\limits_{k=n}^{\infty} \left\vert \Delta_{rr}a_n\right\vert < \varepsilon \qquad mn \ln m \ln n \left\vert c_{mn} \right\vert < \varepsilon, \\
		& \text{ }\underset{j\geq m}{sup} \text{ } j \ln j\sum\limits_{k=n}^{\infty}  \left\vert c_{jk} \right\vert < \varepsilon,
		\qquad \text{ }\underset{k\geq n}{sup} \text{ } k \ln k\sum\limits_{j=m}^{\infty}  \left\vert c_{jk} \right\vert < \varepsilon.
	\end{split}
	\end{equation}
	The above inequalities follows from Lemma \ref{zbieznosc_pojedynczych_z_supremum}, Lemma \ref{zbieznosc_podwojny_bez_supremum},   \eqref{dodatkowy_warunek_dla_max} \eqref{dodatkowy_warunek_dla_max_k}, \eqref{prof 13}.
	
	When $x = 0$ or $y=0$ and $x = \pi$ or $y = \pi$ the proof is trivial.
	Using $\eqref{zbieznosc_w_punktach}$ we obtain that $\eqref{zbieznosc_reszty_szeregu_podwojnego}$ is true if $x = \frac{2l_1\pi}{r}$ and $y=\frac{2l_2\pi}{r}$, where $r \geq 3$ and $l_1,l_2 \in \{  1,...,[\frac{r}{2}]-1 \}$ when $r$ is an even number, and $l_1,l_2 \in \{  1,...,[\frac{r}{2}] \}$ when $r$ is an odd number, respectively. 
	
	Suppose $x \in \left(\frac{2l_1\pi}{r},\frac{(2l_1+1)\pi}{r}\right]$ and $y \in \left(\frac{2l_2\pi}{r},\frac{(2l_2+1)\pi}{r}\right]$, where $l_1,l_2 \in \{  0,1,...,[\frac{r}{2}]-1 \}$ when $r$ is an even number, and $l_1,l_2\in \{ 0,1,...,[\frac{r}{2}] \}$ when $r$ is an odd number, respectively. Let $\mu = \left[\frac{1}{x-\frac{2l_1\pi}{r}}\right]^p$ and $\nu = \left[\frac{1}{y-\frac{2l_2\pi}{r}}\right]^p$. Using elementary calculations we obtain:
	\begin{align}\label{oszacowanie_do_twierdzenia_glownego}
		\sin(jx)&\sin(ky)\notag  \\
		& 
		= \left( \sin(jx) - \sin(j\frac{2l_1\pi}{r}) \right) \left( \sin(ky) - \sin(k\frac{2l_2\pi}{r}) \right)\notag  \\
		&
		+ \left( \sin(jx) - \sin(j\frac{2l_1\pi}{r}) \right) \sin(k\frac{2l_2\pi}{r}) \notag  \\
		&
		+ \sin(j\frac{2l_1\pi}{r})\left( \sin(ky) - \sin(k\frac{2l_2\pi}{r}) \right) + \sin(j \frac{2l_1\pi}{r}) \sin(k \frac{2l_2\pi}{r})
	\end{align}
	Applying Lagrange’s mean value theorem to the functions $f(x) = \sin(jx)$ and $f(y) = \sin(ky)$ on the	interval $\left[\frac{2l_1\pi}{r}, x\right]$ and $\left[\frac{2l_2\pi}{r}, y\right]$ respectively, we obtain that for each $j,k$ there exist $z_j \in \left( \frac{2l_1\pi}{r}, x\right)$ and $ q_k \in \left( \frac{2l_2\pi}{r}, y\right)$ such that:
	\begin{align}
		\label{twierdzenie_lagrangea_pierwszy_przedzial}
		\begin{split}
			\sin(jx) - \sin(\frac{2l_1\pi}{r}) = j \cos(jz_j) (x - \frac{2l_1\pi}{r}),
			\\
			\sin(ky) - \sin(\frac{2l_2\pi}{r}) = k \cos(kq_k) (y - \frac{2l_2\pi}{r}).
		\end{split}
	\end{align}
	We have to consider the following cases:	

	CASE (a): $\eta <m\leq M \leq \mu$ and $\eta <n\leq N \leq \nu$. In CASE (a) we have four subcases to consider:
	\allowdisplaybreaks
	
	Subcase ($a_1$): $m+\mu^{\frac{1}{p}} > M$ and $n+\nu^{\frac{1}{p}} > N$. Using \eqref{zbieznosc_w_punktach}, \eqref{lematy_zbieznosc_w_twierdzeniu} and \eqref{twierdzenie_lagrangea_pierwszy_przedzial} we obtain:
	\begin{flalign*}
	\allowdisplaybreaks
	\left\vert \sum\limits_{j=m}^{M}\sum\limits_{k=n}^{N} \right. & c_{jk} \sin(jx)\sin(ky) \Bigg\vert & \\
	&
	\leq \left\vert  \sum\limits_{j=m}^{M}\sum\limits_{k=n}^{N} \right.  c_{jk}
	\left( \sin(jx) - \sin(j\frac{2l_1\pi}{r}) \right) \left( \sin(ky) - \sin(k\frac{2l_2\pi}{r}) \right) \Bigg\vert \\
	&
	+ \left\vert  \sum\limits_{j=m}^{M}\sum\limits_{k=n}^{N} \right.  c_{jk} \left( \sin(jx) - \sin(j\frac{2l_1\pi}{r}) \right) \sin(k\frac{2l_2\pi}{r}) \Bigg\vert \\
	&
	+ \left\vert  \sum\limits_{j=m}^{M}\sum\limits_{k=n}^{N} \right.  c_{jk} \sin(j\frac{2l_1\pi}{r})\left( \sin(ky) - \sin(k\frac{2l_2\pi}{r}) \right) \Bigg\vert \\
	&
	+ \left\vert  \sum\limits_{j=m}^{M}\sum\limits_{k=n}^{N} \right.  c_{jk} \sin(j \frac{2l_1\pi}{r}) \sin(k \frac{2l_2\pi}{r}) \Bigg\vert \\
	&
	<
	\left(x - \frac{2l_1\pi}{r}\right)\left(y - \frac{2l_2\pi}{r}\right) \sum\limits_{j=m}^{M}\sum\limits_{k=n}^{N} \left\vert j \cos(jz_j)  k \cos(kq_k)  c_{jk}  \right\vert  \\
	& 
	+ \left(x - \frac{2l_1\pi}{r}\right) \sum\limits_{j=m}^{M}\sum\limits_{k=n}^{N}\left\vert c_{jk} j \cos(jz_j) \sin(k\frac{2l_2\pi}{r}) \right\vert \\ 
	&
	+ \left(y - \frac{2l_2\pi}{r}\right) \sum\limits_{j=m}^{M}\sum\limits_{k=n}^{N}\left\vert c_{jk} k \cos(kq_k) \sin(j\frac{2l_1\pi}{r}) \right\vert 
	+ \varepsilon  \\
	&
	\leq \frac{1}{\mu^{\frac{1}{p}}\nu^{\frac{1}{p}}} \sum\limits_{j=m}^{m+\mu^{\frac{1}{p}}-1}\sum\limits_{k=n}^{n+\nu^{\frac{1}{p}}-1} jk\ln j\ln k \left\vert c_{jk}  \right\vert \\
	&
	+ \frac{1}{\mu^{\frac{1}{p}}} \sum\limits_{j=m}^{m+\mu^{\frac{1}{p}}-1} \underset{j\geq m}{sup} \text{ } j \ln j\sum\limits_{k=n}^{\infty}  \left\vert c_{jk}  \right\vert 
	+ 
	\frac{1}{\nu^{\frac{1}{p}}} \sum\limits_{k=n}^{n+\nu^{\frac{1}{p}}-1} \underset{k\geq n}{sup } \text{ } k \ln k\sum\limits_{j=m}^{\infty}  \left\vert c_{jk}  \right\vert	+ \varepsilon < 3 \varepsilon.
	\end{flalign*}

	Subcase ($a_2$): $m+\mu^{\frac{1}{p}} > M$ and $n+\nu^{\frac{1}{p}}\leq N$.	We have:
	\begin{flalign*}
		\sum\limits_{j=m}^{M}\sum\limits_{k=n}^{N} & c_{jk} \sin(jx)\sin(ky)
		= \sum\limits_{j=m}^{M}\sum\limits_{k=n}^{n+\nu^{\frac{1}{p}}-1} c_{jk} \sin(jx)\sin(ky)
		& \\
		& + 			 \sum\limits_{j=m}^{M}\sum\limits_{k=n+\nu^{\frac{1}{p}}}^{N} c_{jk} \sin(jx)\sin(ky)
		=  S_1 + S_2.
	\end{flalign*}
	Similarly as subcase $(a_1)$ we can show that:
	\begin{flalign*}
	\left\vert S_1\right\vert  &=  \left\vert \sum\limits_{j=m}^{M}\sum\limits_{k=n}^{n+\nu^{\frac{1}{p}}-1} c_{jk} \sin(jx)\sin(ky) \right\vert
	\leq \frac{1}{\mu^{\frac{1}{p}}\nu^{\frac{1}{p}}} \sum\limits_{j=m}^{m+\mu^{\frac{1}{p}}-1}\sum\limits_{k=n}^{n+\nu^{\frac{1}{p}}-1} j\ln jk\ln k \left\vert c_{jk}  \right\vert & \\
	&
	+ \frac{1}{\mu^{\frac{1}{p}}} \sum\limits_{j=m}^{m+\mu^{\frac{1}{p}}-1} \underset{j\geq m}{sup} \text{ } j \ln j\sum\limits_{k=n}^{\infty}  \left\vert c_{jk}  \right\vert 
	+ 
	\frac{1}{\nu^{\frac{1}{p}}} \sum\limits_{k=n}^{n+\nu^{\frac{1}{p}}-1} \underset{k\geq n}{sup } \text{ } k \ln k\sum\limits_{j=m}^{\infty}  \left\vert c_{jk}  \right\vert	+ \varepsilon < 3 \varepsilon.
	\end{flalign*}
	Using Lemma \ref{lemat_do_warunku_z_logarytmem}, \eqref{lematy_zbieznosc_w_twierdzeniu}, \eqref{oszacowanie_do_twierdzenia_glownego} and \eqref{twierdzenie_lagrangea_pierwszy_przedzial} we have:
	\allowdisplaybreaks
	\begin{flalign*}
	\left\vert S_2\right\vert & =  \left\vert \sum\limits_{j=m}^{M}\sum\limits_{k=n+\nu^{\frac{1}{p}}}^{N} c_{jk} \sin(jx)\sin(ky) \right\vert
	&\\
	&
	= \left\vert  \sum\limits_{j=m}^{M}\sum\limits_{k=n+\nu^{\frac{1}{p}}}^{N} \right.  c_{jk} \left( j\cos(jz_j)\left(x - \frac{2l_1\pi}{r} \right) + \sin(j\frac{2l_1\pi}{r}) \right) \sin(ky) \Bigg\vert  \\
	&
	\leq \left(x - \frac{2l_1\pi}{r}\right)  \sum\limits_{j=m}^{M}\sum\limits_{k=n+\nu^{\frac{1}{p}}}^{N} \left\vert c_{jk} j \cos(jz_j) \right\vert
	+  \sum\limits_{j=m}^{M}\sum\limits_{k=n+\nu^{\frac{1}{p}}}^{N} \left\vert c_{jk}\sin(j\frac{2l_1\pi}{r}) \sin(ky) \right\vert \\
	&
	< \frac{4\varepsilon}{\mu^{\frac{1}{p}}} \sum\limits_{j=m}^{m+\mu^{\frac{1}{p}}-1} \sum\limits_{k=n+\nu^{\frac{1}{p}}}^{n+\nu-1} \frac{ 1}{(k+1)\ln (k+1)} 
	+ 4  \underset{k\geq n}{\text{ } sup } \text{ } k \ln k\sum\limits_{j=m}^{\infty}  \left\vert c_{jk} \right\vert \sum\limits_{k=n+\nu^{\frac{1}{p}}}^{n+\nu-1} \frac{1}{(k+1)\ln (k+1)} \\
	&
	\leq \frac{4\varepsilon}{\mu^{\frac{1}{p}}} \sum\limits_{j=m}^{m+\mu^{\frac{1}{p}}-1} \sum\limits_{k=n+\nu^{\frac{1}{p}}}^{n+\nu-1} \int\limits_{k}^{k+1}\frac{ 1}{x\ln x} dx
	+ 4  \underset{k\geq n}{ \text{ }sup } \text{ } k \ln k\sum\limits_{j=m}^{\infty}  \left\vert c_{jk}  \right\vert  \sum\limits_{k=n+\nu^{\frac{1}{p}}}^{n+\nu- 1} \int\limits_{k}^{k+1}\frac{ 1}{x\ln x} dx \\
	&
	\leq \frac{4\varepsilon}{\mu^{\frac{1}{p}}} \sum\limits_{j=m}^{m+\mu^{\frac{1}{p}}-1} \int\limits_{n+\nu^{\frac{1}{p}}}^{n+\nu} \frac{ 1}{k\ln k} dk
	+ 4 \varepsilon \int\limits_{n+\nu^{\frac{1}{p}}}^{n+\nu} \frac{ 1}{k\ln k} dk
	\leq \frac{4\varepsilon}{\mu^{\frac{1}{p}}} \sum\limits_{j=m}^{m+\mu^{\frac{1}{p}}-1} \ln p 
	+ 4 \varepsilon \ln p 
	\\
	&
	\leq 8 \varepsilon \ln p.
	\end{flalign*}
	Hence:
	\begin{equation*}
	\left\vert \sum\limits_{j=m}^{M}\sum\limits_{k=n}^{N} c_{jk} \sin(jx)\sin(ky) \right\vert
	<
	\varepsilon \left(3+ 8 \ln p\right).
	\end{equation*}

	Subcase ($a_3$): $m+\mu^{\frac{1}{p}} \leq M$ and $n+\nu^{\frac{1}{p}} > N$. Similarly as in the subcase ($a_2$), we can show:
	\begin{equation*}
	\left\vert \sum\limits_{j=m}^{M}\sum\limits_{k=n}^{N} c_{jk} \sin(jx)\sin(ky) \right\vert
	<
	\varepsilon \left(3+ 8 \ln p\right).
	\end{equation*}

	Subcase ($a_4$): $m+\mu^{\frac{1}{p}} \leq M$ and $n+\nu^{\frac{1}{p}} \leq N$.	We have:
	\begin{flalign*}
	\sum\limits_{j=m}^{M}\sum\limits_{k=n}^{N} & c_{jk} \sin(jx)\sin(ky) &\\
	&
	= \sum\limits_{j=m}^{m+\mu^{\frac{1}{p}}-1}\sum\limits_{k=n}^{n+\nu^{\frac{1}{p}}-1} c_{jk} \sin(jx)\sin(ky) + \sum\limits_{j=m+\mu^{\frac{1}{p}}}^{M}\sum\limits_{k=n}^{n+\nu^{\frac{1}{p}}-1} c_{jk} \sin(jx)\sin(ky) \\
	& 
	+ \sum\limits_{j=m}^{m+\mu^{\frac{1}{p}}-1}\sum\limits_{k=n+\nu^{\frac{1}{p}}}^{N} c_{jk} \sin(jx)\sin(ky) + \sum\limits_{j=m+\mu^{\frac{1}{p}}}^{M}\sum\limits_{k=n+\nu^{\frac{1}{p}}}^{N} c_{jk} \sin(jx)\sin(ky) \\
	&
	= S_3 + S_4 +S_5 + S_6.
	\end{flalign*}
	Similarly as subcase $(a_1)$ we can show that:
	\begin{flalign*}
	\left\vert S_3\right\vert & =  \left\vert \sum\limits_{j=m}^{m+\mu^{\frac{1}{p}}-1}\sum\limits_{k=n}^{n+\nu^{\frac{1}{p}}-1} c_{jk} \sin(jx)\sin(ky) \right\vert
	< 3 \varepsilon .
	\end{flalign*}
	Similarly as in subcase $(a_2)$ we can show that:
	\begin{flalign*}
	\left\vert S_4\right\vert & =  \left\vert \sum\limits_{j=m+\mu^{\frac{1}{p}}}^{M}\sum\limits_{k=n}^{n+\nu^{\frac{1}{p}}-1} c_{jk} \sin(jx)\sin(ky) \right\vert &\\
	&
	= 
	\left\vert  \left( y - \frac{2l_2\pi}{r} \right)\sum\limits_{j=m+\mu^{\frac{1}{p}}}^{M}\sum\limits_{k=n}^{n+\nu^{\frac{1}{p}}-1} c_{jk} \sin(jx)k\cos(kq_k)  + \sum\limits_{j=m+\mu^{\frac{1}{p}}}^{M}\sum\limits_{k=n}^{n+\nu^{\frac{1}{p}}-1} c_{jk} \sin(jx)\sin(k\frac{2l_2\pi}{r})\right\vert \\
	&
	<
	\frac{4\varepsilon}{\nu^{\frac{1}{p}}} \sum\limits_{j=m+\mu^{\frac{1}{p}}}^{m+\mu-1}\sum\limits_{k=n}^{n+\nu^{\frac{1}{p}}-1} \frac{1}{(j+1)\ln(j+1)}
	+ 4 \underset{j\geq m}{\text{ }sup } \text{ } j \ln j \sum\limits_{k=n}^{\infty} \left\vert c_{jk} \right\vert \sum\limits_{j=m+\mu^{\frac{1}{p}}}^{m+\mu-1} \frac{1}{(j+1)\ln(j+1)} 	 \\
	&
	\leq
	8\varepsilon \sum\limits_{j=m+\mu^{\frac{1}{p}}}^{m+\mu-1}	\int\limits_{j}^{j+1}\frac{1}{x \ln x} dx 
	\leq 8 \varepsilon \int\limits_{j+\mu^{\frac{1}{p}}}^{j+\mu}\frac{1}{x \ln x} dx
	\leq 8 \varepsilon \ln p
	\end{flalign*}
	and 
	\begin{equation*}
	\left\vert S_5\right\vert  = \left\vert  \sum\limits_{j=m}^{m+\mu^{\frac{1}{p}}-1}\sum\limits_{k=n+\nu^{\frac{1}{p}}}^{N} c_{jk} \sin(jx)\sin(ky) \right\vert < \varepsilon \left(8 \ln p\right).
	\end{equation*}
	Next, using Lemma \ref{lemat_do_warunku_z_logarytmem} we have:
	\begin{flalign*}
	\left\vert S_6\right\vert & =  \left\vert \sum\limits_{j=m+\mu^{\frac{1}{p}}}^{M}\sum\limits_{k=n+\nu^{\frac{1}{p}}}^{N} c_{jk} \sin(jx)\sin(ky) \right\vert
	\leq \sum\limits_{j=m+\mu^{\frac{1}{p}}}^{M}\sum\limits_{k=n+\nu^{\frac{1}{p}}}^{N} \left\vert c_{jk} \right\vert & \\
	&
	< 16 \varepsilon \sum\limits_{j=m+\mu^{\frac{1}{p}}}^{m+\mu-1} \sum\limits_{k=n+\nu^{\frac{1}{p}}}^{n+\nu-1} \frac{1}{(j+1)\ln (j+1)}  \frac{1}{(k+1)\ln (k+1)} 
	& \\
	&
	\leq 16 \varepsilon \sum\limits_{j=m+\mu^{\frac{1}{p}}}^{m+\mu-1}\int\limits_{j}^{j+1}\frac{1}{x \ln x} dx \sum\limits_{k=n+\nu^{\frac{1}{p}}}^{n+\nu-1} \int\limits_{k}^{k+1}\frac{1}{y \ln y} dy
	& \\
	&
	\leq 16 \varepsilon \int\limits_{n+\nu^{\frac{1}{p}}}^{n+\nu} \int\limits_{m+\mu^{\frac{1}{p}}}^{m+\mu}\frac{1}{x\ln x} \frac{1}{y\ln y} dx dy
	\leq 16 \varepsilon \ln^2p.
	\end{flalign*}	
	
	CASE (b): $max\{\mu,\eta \} <m\leq M$ and $\eta <n\leq N \leq \nu$. In CASE (b) we have two subcases to consider:	
	
	Subcase ($b_1$): $n+\nu^{\frac{1}{p}} > N$. From \eqref{oszacowanie_do_twierdzenia_glownego} and \eqref{twierdzenie_lagrangea_pierwszy_przedzial} we have: 
	\begin{flalign*}
		\sum\limits_{j=m}^{M}\sum\limits_{k=n}^{N} & c_{jk} \sin(jx)\sin(ky) 
		= \sum\limits_{k=n}^{N}  \left( k \cos(kq_k) \left( y - \frac{2l_2\pi}{r} \right) + \sin(k\frac{2l_2\pi}{r}) \right)\sum\limits_{j=m}^{M} c_{jk}\sin(jx) & \\
		&
		= \left( y - \frac{2l_2\pi}{r} \right) \sum\limits_{k=n}^{N} k\cos(kq_k) \sum\limits_{j=m}^{M} c_{jk}\sin(jx) 
		+ \sum\limits_{k=n}^{N} \sum\limits_{j=m}^{M} c_{jk}\sin(k\frac{2l_2\pi}{r})\sin(jx)& \\
		&
		= S_7 + S_8.	
	\end{flalign*}
	From Lemma \ref{lemat_do_glownego_tw} and  \eqref{lematy_zbieznosc_w_twierdzeniu} we have:	
	\begin{flalign}\label{s_7}
		\left\vert S_{7} \right\vert &
		\leq \left( y - \frac{2l_2\pi}{r} \right) \sum\limits_{k=n}^{n+\nu^{\frac{1}{p}}-1} k \left\vert \sum\limits_{j=m}^{M} c_{jk}\sin(jx) \right\vert \notag & \\
		&
		\leq \frac{1}{\nu^{\frac{1}{p}}}
		\sum\limits_{k=n}^{n+\nu^{\frac{1}{p}}-1} \frac{k}{2\left(\frac{rx}{\pi}-2l_1\right)} \left( \sum\limits_{j=m}^{M} \left\vert\Delta_{r0}c_{jk}\right\vert +
		\sum\limits_{j=M+1}^{M+r} \left\vert c_{jk}\right\vert +
		\sum\limits_{j=m}^{m+r-1}\left\vert c_{jk}\right\vert \right) \notag & \\
		&
		\leq
		\frac{\pi (\mu +1)^{\frac{1}{p}}}{2 r \nu^{\frac{1}{p}} } \sum\limits_{k=n}^{n+\nu^{\frac{1}{p}}-1} k \left( \sum\limits_{j=m}^{M} \left\vert\Delta_{r0}c_{jk}\right\vert +
		\sum\limits_{j=M+1}^{M+r} \left\vert c_{jk}\right\vert +
		\sum\limits_{j=m}^{m+r-1}\left\vert c_{jk}\right\vert \right) \notag \\
		&
		\leq
		\frac{\pi}{2\nu^{\frac{1}{p}} r} \sum\limits_{k=n}^{n+\nu^{\frac{1}{p}}-1} \left(m^{\frac{1}{p}} \text{ }\underset{k\geq n}{sup} \text{ }k  \sum\limits_{j=m}^{M} \left\vert\Delta_{r0}c_{jk}\right\vert +
		\sum\limits_{j=M+1}^{M+r} jk\ln j\ln k \left\vert c_{jk}\right\vert +
		\sum\limits_{j=m}^{m+r-1}jk\ln j\ln k \left\vert c_{jk}\right\vert \right) \notag  \\
		&
		<
		\frac{\pi}{2\nu^{\frac{1}{p}} r} \sum\limits_{k=n}^{n+\nu^{\frac{1}{p}}-1} \left( \varepsilon + r\varepsilon + r\varepsilon \right) 
		\leq 3 \pi\varepsilon.
	\end{flalign}
	Next, using Lemma \ref{lemat_do_glownego_tw}, \eqref{lematy_zbieznosc_w_twierdzeniu} and the definition of the class $DGM\left(p,{}_2\alpha, {}_2\beta,{}_2\gamma,r\right) $ we have:	
	\begin{flalign}\label{s_8}
		\left\vert S_{8} \right\vert & \leq \sum\limits_{k=n}^{n+\nu^{\frac{1}{p}}-1} \left\vert \sum\limits_{j=m}^{M} c_{jk}\sin(jx) \right\vert &\notag \\ 
		&
		\leq \sum\limits_{k=n}^{n+\nu^{\frac{1}{p}}-1} \frac{1}{2\left(\frac{rx}{\pi}-2l_1\right)} \left( \sum\limits_{j=m}^{M} \left\vert\Delta_{r0}c_{jk}\right\vert +
		\sum\limits_{j=M+1}^{M+r} \left\vert c_{jk}\right\vert +
		\sum\limits_{j=m}^{m+r-1}\left\vert c_{jk}\right\vert \right) \notag \\ 
		&
		\leq
		\frac{\pi (\mu +1)^{\frac{1}{p}}}{2 r} \sum\limits_{k=n}^{n+\nu^{\frac{1}{p}}-1} \left( \sum\limits_{j=m}^{M} \left\vert\Delta_{r0}c_{jk}\right\vert +
		\sum\limits_{j=M+1}^{M+r} \left\vert c_{jk}\right\vert +
		\sum\limits_{j=m}^{m+r-1}\left\vert c_{jk}\right\vert \right)\notag \\
		&
		\leq
		\frac{\pi}{ 2r} \left( m^{\frac{1}{p}}  \sum\limits_{k=n}^{\infty} \sum\limits_{r=0}^{\infty} \sum\limits_{j=2^{r} m}^{2^{r+1} m} \left\vert\Delta_{r0}c_{jk}\right\vert 
		+  \sum\limits_{j=M+1}^{M+r} \text{ }\underset{j\geq m}{sup} \text{ }  j\ln j \sum\limits_{k=n}^{\infty} \left\vert c_{jk}\right\vert \right. \notag \\
		&
		\left. +
		\sum\limits_{j=m}^{m+r-1} \text{ }\underset{j\geq m}{sup} \text{ } j\ln j \sum\limits_{k=n}^{\infty} \left\vert c_{jk}\right\vert \right) \notag \\
		&
		< 
		\frac{\pi}{ 2r} \left(m^{\frac{1}{p}}  \sum\limits_{k=n}^{\infty} \sum\limits_{r=0}^{\infty} \left(2^{r}m\right)^{1-\frac{1}{p}} \left( \sum\limits_{j=2^{r} m}^{2^{r+1} m-1} \left\vert\Delta_{r0}c_{jk}\right\vert^p \right)^\frac{1}{p}
		+ 2 r \varepsilon \right) \notag \\
		&
		\leq 
		\frac{\pi}{2r} \left(m^{\frac{1}{p}}  \sum\limits_{k=n}^{\infty} \sum\limits_{r=0}^{\infty}  
		\frac{C}{\left( 2^{r}m\right)^{\frac{1}{p}}} \text{ }\underset{b_1(2^{r}m)\leq M \leq \lambda b_1(2^{r}m)}{max} \text{ } \sum\limits_{j=M}^{2M} \left\vert c_{jk}\right\vert 
		+ 2 r \varepsilon \right) \notag \\
		&
		\leq 
		\frac{\pi}{ 2r} \left( C  \sum\limits_{r=0}^{\infty}  
		\frac{1}{2^{\frac{r}{p}}} \sum\limits_{j=b_1(2^{r}m)}^{2\lambda b_1(2^{r}m)} \frac{1}{j\ln j} \left( \text{ }\underset{j\geq m}{sup} \text{ } j\ln j\sum\limits_{k=n}^{\infty} \left\vert c_{jk}\right\vert  \right)
		+ 2 r \varepsilon \right)\notag \\
		&
		\leq 
		\frac{\pi\varepsilon}{ 2r} \left( C  \sum\limits_{r=0}^{\infty}  
		\frac{1}{2^{\frac{r}{p}}} \sum\limits_{j=b_1(2^{r}m)}^{2\lambda b_1(2^{r}m)} \frac{1}{j\ln j}
		+ 2 r \varepsilon \right)
		\leq 
		\frac{\pi\varepsilon}{r} \left( C\lambda  \sum\limits_{r=0}^{\infty}  
		\frac{1}{2^\frac{r}{p}}	+ r \right)
		\ll \varepsilon
	\end{flalign}
	Hence \eqref{zbieznosc_reszty_szeregu_podwojnego} holds.

	Subcase ($b_2$): $n+\nu^{\frac{1}{p}} \leq N$. We have:
	
	\begin{flalign*}
	\sum\limits_{j=m}^{M}\sum\limits_{k=n}^{N} &  c_{jk} \sin(jx)\sin(ky) 
	= \sum\limits_{k=n}^{n+\nu^{\frac{1}{p}}-1} \sin(ky) \sum\limits_{j=m}^{M} c_{jk}\sin(jx) + \sum\limits_{k=n+\nu^{\frac{1}{p}}}^{N} \sin(ky) \sum\limits_{j=m}^{M} c_{jk}\sin(jx) \notag & \\
	&
	= \sum\limits_{k=n}^{n+\nu^{\frac{1}{p}}-1}  \left( \left( k \cos(kq_k) \left( y - \frac{2l_2\pi}{r} \right)\right) + \sin(k\frac{2l_2\pi}{r}) \right)\sum\limits_{j=m}^{M} c_{jk}\sin(jx) \notag &\\
	&
	+ \sum\limits_{k=n+\nu^{\frac{1}{p}}}^{N} \sin(ky) \sum\limits_{j=m}^{M} c_{jk}\sin(jx) \notag & \\
	&
	= \left( y - \frac{2l_2\pi}{r} \right) \sum\limits_{k=n}^{n+\nu^{\frac{1}{p}}-1} k\cos(kq_k) \sum\limits_{j=m}^{M} c_{jk}\sin(jx) 
	+ \sum\limits_{k=n}^{n+\nu^{\frac{1}{p}}-1} \sum\limits_{j=m}^{M} c_{jk}\sin(k\frac{2l_2\pi}{r})\sin(jx) \notag &\\
	&
	+ \sum\limits_{k=n+\nu^{\frac{1}{p}}}^{N} \sin(ky) \sum\limits_{j=m}^{M} c_{jk}\sin(jx) = S_9 + S_{10} + S_{11}.
	\end{flalign*}
	Similarly as in \eqref{s_7} and \eqref{s_8} we can show that $\left\vert S_9 \right\vert \leq 3 \pi \varepsilon$ and $\left\vert S_{10} \right\vert \ll \varepsilon$.
	Using Lemma \ref{lemat_do_warunku_z_logarytmem}, Lemma \ref{zbieznosc_pojedynczych_z_supremum}, \eqref{lematy_zbieznosc_w_twierdzeniu}, the definition of the class $DGM\left(p,{}_2\alpha, {}_2\beta,{}_2\gamma,r\right) $ and H\"{o}lder's inequality with $p>1$ we obtain:
	\begin{flalign*}
	\left\vert S_{11} \right\vert & \leq \sum\limits_{k=n+\nu^{\frac{1}{p}}}^{N} \left\vert\sum\limits_{j=m}^{M} c_{jk}\sin(jx)  \right\vert & \\
	&
	\leq
	\sum\limits_{k=n+\nu^{\frac{1}{p}}}^{N} \frac{1}{2\left(\frac{rx}{\pi}-2l_1\right)} \left( \sum\limits_{j=m}^{M} \left\vert\Delta_{r0}c_{jk}\right\vert  +
	\sum\limits_{j=M+1}^{M+r} \left\vert c_{jk}\right\vert +
	\sum\limits_{j=m}^{m+r-1}\left\vert c_{jk}\right\vert \right) \\
	&
	\leq
	\frac{\pi (\mu +1)^{\frac{1}{p}}}{ 2r}\sum\limits_{k=n+\nu^{\frac{1}{p}}}^{N} \left( \sum\limits_{j=m}^{M} \left\vert\Delta_{r0}c_{jk}\right\vert +
	\sum\limits_{j=M+1}^{M+r} \left\vert c_{jk}\right\vert +
	\sum\limits_{j=m}^{m+r-1}\left\vert c_{jk}\right\vert \right) \\
	&
	< 
	\frac{\pi}{ 2r} m^{\frac{1}{p}}  \sum\limits_{k=n}^{\infty} \sum\limits_{r=0}^{\infty} \left(2^{r}m\right)^{1-\frac{1}{p}} \left( \sum\limits_{j=2^{r} m}^{2^{r+1} m-1} \left\vert\Delta_{r0}c_{jk}\right\vert^p \right)^\frac{1}{p}
	+ \frac{8 \pi r \varepsilon}{2r}  \sum\limits_{k=n+\nu^{\frac{1}{p}}}^{n+\nu-1}  \frac{1}{(k+1) \ln (k+1)} \\
	&
	\leq 
	\frac{\pi}{2r} m^{\frac{1}{p}}  \sum\limits_{k=n}^{\infty} \sum\limits_{r=0}^{\infty}  
	\frac{C}{\left( 2^{r}m\right)^{\frac{1}{p}}} \text{ }\underset{b_1(2^{r}m)\leq M \leq \lambda b_1(2^{r}m)}{max} \text{ } \sum\limits_{j=M}^{2M} \left\vert c_{jk}\right\vert 
	+ 4 \pi \varepsilon  \sum\limits_{k=n+\nu^{\frac{1}{p}}}^{n+\nu-1} \int\limits_{k}^{k+1} \frac{1}{x \ln x} dx   \\
	&
	\leq 
	\frac{C\pi}{ 2r}  \sum\limits_{r=0}^{\infty}  
	\frac{1}{2^{\frac{r}{p}}} \sum\limits_{j=b_1(2^{r}m)}^{2\lambda b_1(2^{r+1}m)} \frac{1}{j \ln j} \left( \text{ }\underset{j\geq m}{sup} \text{ } j\ln j \sum\limits_{k=n}^{\infty} \left\vert c_{jk}\right\vert  \right)
	+ 4 \pi \varepsilon  \int\limits_{n+\nu^{\frac{1}{p}}}^{n+\nu} \frac{1}{x \ln x} dx \\
	&
	\leq 
	\frac{C\pi\varepsilon}{ 2r}  \sum\limits_{r=0}^{\infty}  
	\frac{1}{2^{\frac{r}{p}}} \sum\limits_{j=b_1(2^{r}m)}^{2\lambda b_1(2^{r+1}m)} \frac{1}{j \ln j}
	+ 4 \pi \varepsilon \ln p  \\
	&
	\leq 
	\frac{C\pi\varepsilon2\lambda}{ 2r}  \sum\limits_{r=0}^{\infty}  
	\frac{1}{2^\frac{r}{p}}	+ 4 \pi \varepsilon \ln p
	\ll \varepsilon.
	\end{flalign*}
	So we have that \eqref{zbieznosc_reszty_szeregu_podwojnego} is satisfied.

	CASE (c):  $\eta <m\leq M \leq \mu$ and $max\{\eta,\nu\} <n\leq N$. This is the symmetric counterpart of CASE (b), and the proof is similar. Therefore \eqref{zbieznosc_reszty_szeregu_podwojnego} holds.
	
	CASE (d): $max\{\eta,\mu\} <m\leq M$ and $max\{\eta,\nu\} <n\leq N$. Using Lemma \ref{roznicowanie} we get:
	\begin{flalign*}
	\left\vert \sum\limits_{j=m}^{M}\sum\limits_{k=n}^{N} c_{jk} \right. & \sin(jx)\sin(ky) \bigg\vert = 
	\left\vert \sum\limits_{j=m}^{M} \sin(jx) \sum\limits_{k=n}^{N} c_{jk} \sin(ky) \right\vert & \\
	&
	=
	\left\vert \sum\limits_{j=m}^{M} \sin(jx) \left( - \sum\limits_{k=n}^{N} \Delta_{0r} c_{jk} \tilde{D}_{k,r}(y)
	+
	\sum\limits_{k=N+1}^{N+r} \left\vert c_{jk} \right\vert \left\vert \tilde{D}_{k,-r}(y)\right\vert \right. \right. \\
	&
	\left. \left.
	-
	\sum\limits_{k=n}^{n+r-1} c_{jk}  \tilde{D}_{k,-r}(y) \right) \right\vert \\
	&
	\leq
	\left\vert \sum\limits_{j=m}^{M}\sum\limits_{k=n}^{N}  \Delta_{rr} c_{jk} \tilde{D}_{j,r}(x) \tilde{D}_{k,r}(y)\right\vert
	+
	\left\vert \sum\limits_{j=M+1}^{M+r}\sum\limits_{k=n}^{N}  \Delta_{0r} c_{jk} \tilde{D}_{j,-r}(x) \tilde{D}_{k,r}(y)\right\vert \\
	&
	+
	\left\vert \sum\limits_{j=m}^{m+r-1}\sum\limits_{k=n}^{N}  \Delta_{0r} c_{jk} \tilde{D}_{j,-r}(x) \tilde{D}_{k,r}(y)\right\vert
	+
	\left\vert \sum\limits_{j=m}^{M}\sum\limits_{k=N+1}^{N+r}  \Delta_{r0} c_{jk} \tilde{D}_{j,r}(x) \tilde{D}_{k,-r}(y)\right\vert \\
	&
	+
	\left\vert \sum\limits_{j=m}^{M}\sum\limits_{k=n}^{n+r-1} \Delta_{r0} c_{jk} \tilde{D}_{j,r}(x) \tilde{D}_{k,-r}(y)\right\vert
	+
	\left\vert \sum\limits_{j=M+1}^{M+r}\sum\limits_{k=N+1}^{N+r} c_{jk} \tilde{D}_{j,-r}(x) \tilde{D}_{k,-r}(y)\right\vert \\
	&
	+
	\left\vert \sum\limits_{j=m}^{m+r-1}\sum\limits_{k=N+1}^{N+r} c_{jk} \tilde{D}_{j,-r}(x) \tilde{D}_{k,-r}(y)\right\vert
	+
	\left\vert \sum\limits_{j=M}^{M+r}\sum\limits_{k=n}^{n+r-1} c_{jk} \tilde{D}_{j,-r}(x) \tilde{D}_{k,-r}(y)\right\vert \\
	&
	+
	\left\vert \sum\limits_{j=m}^{m+r-1}\sum\limits_{k=n}^{n+r-1} c_{jk} \tilde{D}_{j,-r}(x) \tilde{D}_{k,-r}(y)\right\vert.
	\end{flalign*}
	Using \eqref{oszacowanie_jadro_diri} and \eqref{lematy_zbieznosc_w_twierdzeniu} we have:
	\begin{flalign*}
	\left\vert \sum\limits_{j=m}^{M}\sum\limits_{k=n}^{N} \right. & c_{jk} \sin(jx)\sin(ky) \bigg\vert 
	\leq \frac{\pi^2}{r^2} (\mu+1)^\frac{1}{p}(\nu+1)^\frac{1}{p} \left( \sum\limits_{j=m}^{M}\sum\limits_{k=n}^{N}  \left\vert \Delta_{rr} c_{jk}\right\vert 
	+ \sum\limits_{j=M+1}^{M+r}\sum\limits_{k=n}^{N} \left\vert \Delta_{0r} c_{jk}\right\vert
	\right. & \\
	&
	+ \sum\limits_{j=m}^{m+r-1}\sum\limits_{k=n}^{N} \left\vert \Delta_{0r} c_{jk}\right\vert
	+ \sum\limits_{j=m}^{M}\sum\limits_{k=N+1}^{N+r} \left\vert \Delta_{r0} c_{jk}\right\vert
	+ \sum\limits_{j=m}^{M}\sum\limits_{k=n}^{n+r-1} \left\vert \Delta_{r0} c_{jk}\right\vert
	+ \sum\limits_{j=M+1}^{M+r}\sum\limits_{k=N+1}^{N+r} \left\vert c_{jk}\right\vert \\
	&
	\left.
	+ \sum\limits_{j=m}^{m+r-1}\sum\limits_{k=N+1}^{N+r} \left\vert c_{jk}\right\vert
	+ \sum\limits_{j=M}^{M+r}\sum\limits_{k=n}^{n+r-1} \left\vert c_{jk}\right\vert
	+ \sum\limits_{j=m}^{m+r-1}\sum\limits_{k=n}^{n+r-1} \left\vert c_{jk}\right\vert \right) \\
	&
	\leq \frac{\pi^2}{r^2}\left( m^\frac{1}{p}n^\frac{1}{p} \sum\limits_{j=m}^{\infty}\sum\limits_{k=n}^{\infty}  \left\vert \Delta_{rr} c_{jk}\right\vert
	+ 2rm^\frac{1}{p} \text{ }\underset{k\geq n}{sup} \text{ }k \sum\limits_{j=m}^{\infty} \left\vert \Delta_{r0} c_{jk}\right\vert
	\right. \\
	&
	\left.
	+2rn^\frac{1}{p}\text{ }\underset{j\geq m}{sup} \text{ }j \sum\limits_{k=n}^{\infty} \left\vert \Delta_{0r} c_{jk}\right\vert
	+ 4r^2 \text{ }\underset{j\geq m \text{ }k\geq n}{sup} \text{ }jk\ln j\ln k \left\vert  c_{jk}\right\vert
	\right) \\
	&
	\leq \frac{\pi^2}{r^2} \left( \varepsilon + 2r \varepsilon + 2r \varepsilon + 4 r^2 \varepsilon \right) \ll \varepsilon. 
	\end{flalign*}
	
	CASE (e): $\eta <m\leq \mu < M$ and $\eta <n\leq N \leq \nu$. This can be split into two cases which are covered by Cases (a) and (b), respectively. Thus \eqref{zbieznosc_reszty_szeregu_podwojnego} holds.
	
	CASE (f): $\eta <m \leq  M \leq \mu$ and $\eta <n\leq\nu < N $. This can be split into two cases which are covered by Cases (a) and (c). Therefore we obtain the same estimate as in Case (e).
	
	CASE (g): $\eta <m \leq \mu < M $ and $max\{\eta,\nu\} <n\leq  N $. This can be split into two cases which are covered by Cases (c) and (d). Thus \eqref{zbieznosc_reszty_szeregu_podwojnego} holds.
	
	CASE (h): $max\{\eta,\mu\} <m \leq  M $ and $\eta <n\leq\nu < N $. This can be split into two cases which are covered by Cases (b) and (d). Therefore we obtain the  estimate \eqref{zbieznosc_reszty_szeregu_podwojnego}.
	
	CASE (i): $\eta <m \leq \mu <  M $ and $\eta <n\leq\nu < N $. This can be split into four cases which are covered by Cases (a) - (d). Therefore \eqref{zbieznosc_reszty_szeregu_podwojnego} is satisfied.
	
	Suppose $x \in (\frac{(2l_1+1)\pi}{r}, \frac{2(l_1+1)\pi}{r})$, $y \in (\frac{2l_2\pi}{r},\frac{(2l_2+1)\pi}{r}]$, $\mu = \left[\frac{1}{\frac{2(l_1+1)\pi }{r}-x}\right]^p$ and $\nu = \left[\frac{1}{y-\frac{2l_2\pi}{r}}\right]^p$, where $l_1 \in \{ 0, 1,...,[\frac{r}{2}]-1 \}$ and $l_2 \in \{ 0, 1,...,[\frac{r}{2}]\}$ when $r \geq 3$ is an odd number, and $l_1,l_2 \in \{  0, 1,...,[\frac{r}{2}] -1 \}$ when $r$ is an even number.
	Using elementary calculations we obtain:
	\begin{align}\label{oszacowanie_do_twierdzenia_glownego_2}
		\sin(jx)&\sin(ky) \notag \\
		& 
		= \left( \sin(jx) - \sin(j\frac{2(l_1+1)\pi}{r}) \right) \left( \sin(ky) - \sin(k\frac{(2l_2+1)\pi}{r}) \right) \notag \\
		&
		+ \left( \sin(jx) - \sin(j\frac{2(l_1+1)\pi}{r}) \right) \sin(k\frac{(2l_2+1)\pi}{r})  \notag\\
		&
		+ \sin(j\frac{2(l_1+1)\pi}{r})\left( \sin(ky) - \sin(k\frac{(2l_2+1)\pi}{r}) \right) \notag \\
		&
		+ \sin(j \frac{2(l_1+1)\pi}{r}) \sin(k \frac{(2l_2+1)\pi}{r})
	\end{align}
	Applying Lagrange’s mean value theorem to the functions $f(x) = \sin(jx)$ and $f(y) = \sin(ky)$ on the	interval $\left[ x, \frac{2(l_1+1)\pi}{r} \right]$ and $\left[\frac{(2l_2+1)\pi}{r}, y\right]$ respectively, we obtain that for each $j,k$ there exist $z_j \in \left( x, \frac{2(l_1+1)\pi}{r} \right),$ and $ q_k \in \left( \frac{(2l_2+1)\pi}{r}, y\right)$ such that
	\begin{align}
		\label{twierdzenie_lagrangea_drugi_przedzial}
		\begin{split}
			\sin(\frac{2(l_1+1)\pi}{r}) - \sin(jx) = j \cos(jz_j) (\frac{2(l_1+1)\pi}{r} - x),
			\\
			\sin(ky) - \sin(\frac{(2l_2+1)\pi}{r}) = k \cos(kq_k) (y - \frac{(2l_2+1)\pi}{r}).
		\end{split}
	\end{align}
	Next, we have to consider nine cases - CASE(a*)-(i*) defined in the same way as CASE (a)-(i). Using \eqref{eq_roznicowanie_druga_polowka}, \eqref{oszacowanie_do_twierdzenia_glownego_2} and \eqref{twierdzenie_lagrangea_drugi_przedzial} instead of \eqref{eq_roznicowanie_pierwsza_polowka}, \eqref{oszacowanie_do_twierdzenia_glownego} and \eqref{twierdzenie_lagrangea_pierwszy_przedzial} respectively, we can show similarly as in Cases (a)-(i) that \eqref{zbieznosc_reszty_szeregu_podwojnego} holds.

	Now, let $x \in (\frac{2l_1\pi}{r},\frac{(2l_1+1)\pi}{r}]$, $y \in (\frac{(2l_2+1)\pi}{r},\frac{2(l_2+1)\pi}{r})$, $\mu = \left[\frac{1}{x-\frac{2l_1\pi}{r}}\right]^p$ and $\nu = \left[\frac{1}{\frac{2(l_2+1)\pi}{r}-y}\right]^p$, where $l_1 \in \{  0, 1,...,[\frac{r}{2}]\}$ and $l_2 \in \{  0, 1,...,[\frac{r}{2}]-1\}$ when $r \geq 3$ is an odd number, and $l_1,l_2 \in \{  0, 1,...,[\frac{r}{2}] -1 \}$ when $r$ is an even number. We have also to consider nine cases - CASE(a**) - (i**). This cases are the symmetric counterparts of cases CASE(a*) - (i*), and the proof is similar. So, we obtain \eqref{zbieznosc_reszty_szeregu_podwojnego}.
	
	Next, let $x \in (\frac{(2l_1+1)\pi}{r},\frac{2(l_1+1)\pi}{r})$ and $y \in (\frac{(2l_2+1)\pi}{r},\frac{2(l_2+1)\pi}{r})$ and $\mu = \left[\frac{1}{\frac{2(l_1+1)\pi}{r}-x}\right]^p$ and $\nu = \left[\frac{1}{\frac{2(l_2+1)\pi}{r}-y}\right]^p$, where $l_1, l_2 \in \{  0, 1,...,[\frac{r}{2}]-1 \}$ for $r \geq 2$. Similarly as above, we can show \eqref{zbieznosc_reszty_szeregu_podwojnego}.
	
	Now, let $y = \frac{2l_2\pi}{r}$, $ x \in  (\frac{2l_1\pi}{r},\frac{(2l_1+1)\pi}{r}]$ and let $\mu = \left[\frac{1}{x-\frac{2l_1\pi}{r}}\right]^p$, where $l_1\in\{0, 1,...,[\frac{r}{2}]-1\}$ and $ l_2 \in \{1,...,[\frac{r}{2}]-1 \}$ when $r \geq 4$ is an even number, and $l_1 \in \{  0, 1,...,[\frac{r}{2}] \}$ and $l_2\in \{ 1,...,[\frac{r}{2}] \}$ when $r \geq 3$ is an odd number. We have to consider the following cases:
	
	CASE (a$^\circ$): $\eta <m\leq M \leq \mu$. In CASE (a$^\circ$) we have to consider two subcase:
	
	Subcase ($a_1^{\circ}$): $m+\mu^{\frac{1}{p}}  > M$.	Using \eqref{dodatkowy_warunek_dla_max} and \eqref{twierdzenie_lagrangea_pierwszy_przedzial} we have:
	\begin{flalign*}
		\left\vert \sum\limits_{j=m}^{M}\sum\limits_{k=n}^{N} \right. & c_{jk} \sin (jx) \sin(k\frac{2l_2\pi}{r})  \bigg\vert &\\
		& 
		= \left\vert \sum\limits_{j=m}^{M}\sum\limits_{k=n}^{N} c_{jk} \left( \left( \sin(jx) - \sin(j\frac{2l_1\pi}{r}) \right) + \sin(j\frac{2l_1\pi}{r}) \right) \sin(k\frac{2l_2\pi}{r})  \right\vert \\
		&
		\leq \left(x - \frac{2l_1 \pi}{r}\right) \sum\limits_{j=m}^{M}\sum\limits_{k=n}^{N} j \left\vert c_{jk} \right\vert 
		+ \left\vert \sum\limits_{j=m}^{M}\sum\limits_{k=n}^{N} c_{jk}\right\vert \\
		&
		< \frac{1}{\mu^{\frac{1}{p}}} \sum\limits_{j=m}^{m+\mu^{\frac{1}{p}}-1} \text{ }\underset{j\geq m}{sup} \text{ } j \ln j \sum\limits_{k=n}^{\infty} \left\vert c_{jk} \right\vert  + \varepsilon
		< 2 \varepsilon.
	\end{flalign*}
	
	Subcase ($a_2^\circ$): $m+\mu^{\frac{1}{p}} \leq M$.	Using \eqref{dodatkowy_warunek_dla_max} and \eqref{twierdzenie_lagrangea_pierwszy_przedzial} we have
	\begin{flalign*}
		\sum\limits_{j=m}^{M}\sum\limits_{k=n}^{N} & c_{jk} \sin (jx) \sin(k\frac{2l_2\pi}{r}) 
		= \sum\limits_{j=m}^{m+\mu^{\frac{1}{p}}-1}\sum\limits_{k=n}^{N} c_{jk} \sin (jx) \sin(k\frac{2l_2\pi}{r})   & \\
		&
		+  \sum\limits_{j=m+\mu^{\frac{1}{p}}}^{M}\sum\limits_{k=n}^{N} c_{jk} \sin (jx) \sin(k\frac{2l_2\pi}{r}) = S_{12} + S_{13}.
	\end{flalign*}
	
	Similarly as in Subcase ($a_1^\circ$), using \eqref{dodatkowy_warunek_dla_max} and \eqref{twierdzenie_lagrangea_pierwszy_przedzial} we can show that:
	\begin{flalign*}
		\left\vert S_{12} \right\vert 
		\ll \frac{1}{\mu^{\frac{1}{p}}} \sum\limits_{j=m}^{m+\mu^{\frac{1}{p}}-1} \text{ }\underset{j\geq m}{sup} \text{ } j \ln j\sum\limits_{k=n}^{\infty} \left\vert c_{jk} \right\vert
		\leq \varepsilon.
	\end{flalign*}
	Using \eqref{dodatkowy_warunek_dla_max} and Lemma \ref{lemat_do_warunku_z_logarytmem}, we get:
	\begin{flalign*}
		\left\vert S_{13} \right\vert &
		\leq  \left\vert \sum\limits_{j=m+\mu^{\frac{1}{p}}}^{m+\mu}\sum\limits_{k=n}^{N} c_{jk} \sin (jx) \sin(k\frac{2l_2\pi}{r})  \right\vert
		\leq  4 \sum\limits_{j=m+\mu^{\frac{1}{p}}}^{m+\mu-1} \frac{1}{(j+1) \ln (j+1)}  \text{ }\underset{j\geq m}{sup} \text{ } j\ln j \sum\limits_{k=n}^{\infty} \left\vert c_{jk} \right\vert & \\
		&
		\leq 4\varepsilon \sum\limits_{j=m+\mu^{\frac{1}{p}}}^{m+\mu-1} \int\limits_{j}^{j+1} \frac{1}{x \ln x} dx
		\leq 4\varepsilon \int\limits_{m+\mu^{\frac{1}{p}}}^{m+\mu} \frac{1}{x \ln x} dx
		\leq 4\varepsilon \ln p.
	\end{flalign*}
	
	CASE (b$^\circ$): $max\{\mu,\eta \} <m\leq M$. Using \eqref{dodatkowy_warunek_dla_max}, \eqref{lematy_zbieznosc_w_twierdzeniu} and the definition of the class $DGM\left(p,{}_2\alpha, {}_2\beta,{}_2\gamma,r\right) $ we have:
	\begin{flalign*}
		\left\vert \sum\limits_{j=m}^{M} \sum\limits_{k=n}^{N} \right. & c_{jk} \sin (jx) \sin(k\frac{2l_2\pi}{r})  \bigg\vert
		\leq \sum\limits_{k=n}^{N} 	\left\vert  \sum\limits_{j=m}^{M} c_{jk} \sin (jx) \right\vert & \\
		&
		\leq \sum\limits_{k=n}^{N} \frac{1}{2\left(\frac{rx}{\pi}-2l_1\right)}  \left( \sum\limits_{j=m}^{M} \left\vert\Delta_{r0}c_{jk}\right\vert + \sum\limits_{j=M+1}^{M+r} \left\vert c_{jk}\right\vert +
		\sum\limits_{j=m}^{m+r-1}\left\vert c_{jk}\right\vert \right) \\
		&
		\leq \frac{\pi (\mu +1)^{\frac{1}{p}}}{ 2r } \sum\limits_{k=n}^{N} \left( \sum\limits_{j=m}^{M} \left\vert\Delta_{r0}c_{jk}\right\vert 
		+
		\sum\limits_{j=M+1}^{M+r} \left\vert c_{jk}\right\vert +
		\sum\limits_{j=m}^{m+r-1}\left\vert c_{jk}\right\vert \right) \\
		&
		\leq \frac{\pi}{ 2r } m^{\frac{1}{p}} \sum\limits_{k=n}^{N} \left( \sum\limits_{j=m}^{\infty} \left\vert\Delta_{r0}c_{jk}\right\vert 
		+ \sum\limits_{j=M+1}^{M+r} \left\vert c_{jk}\right\vert 
		+ \sum\limits_{j=m}^{m+r-1}\left\vert c_{jk}\right\vert \right) \\
		&
		\leq 
		\frac{\pi}{ 2r } m^{\frac{1}{p}}  \sum\limits_{k=n}^{\infty} \sum\limits_{r=0}^{\infty} \left( 2^rm \right)^{1-\frac{1}{p}} \left( \sum\limits_{j=2^{r} m}^{2^{r+1} m} \left\vert\Delta_{r0}c_{jk}\right\vert^p \right)^{\frac{1}{p}} 
		+ \pi \text{ }\underset{j\geq m}{sup} \text{ }  j \ln j \sum\limits_{k=n}^{\infty}  \left\vert c_{jk}\right\vert \\
		&
		\leq 
		\frac{\pi}{ 2r } m^{\frac{1}{p}}  \sum\limits_{k=n}^{\infty} \sum\limits_{r=0}^{\infty}  
		\frac{C}{2^{\frac{r}{p}}} \text{ }\underset{b_1(2^{r}m)\leq M \leq \lambda b_1(2^{r}m)}{max} \text{ } \sum\limits_{j=M}^{2M} \left\vert c_{jk}\right\vert 
		+ \pi \varepsilon  \\
		&
		\leq 
		\frac{C\pi}{ 2r }  \sum\limits_{r=0}^{\infty}  
		\frac{1}{2^{\frac{r}{p}}} \sum\limits_{j=b_1(2^{r}m)}^{2\lambda b_1(2^{r}m)} \frac{1}{j \ln j} \left( \text{ }\underset{j\geq m}{sup} \text{ } j\ln j \sum\limits_{k=n}^{\infty} \left\vert c_{jk}\right\vert  \right)
		+ \pi \varepsilon  \\
		&
		\leq 
		\frac{C\pi\varepsilon}{ 2r}  \sum\limits_{r=0}^{\infty}  
		\frac{1}{2^{\frac{r}{p}}} \sum\limits_{j=b_1(2^{r}m)}^{2\lambda b_1(2^{r}m)} \frac{1}{j \ln j}
		+ \pi \varepsilon  \\
		&
		\leq 
		\frac{C\pi\varepsilon\lambda}{ r}  \sum\limits_{r=0}^{\infty}  
		\frac{1}{2^{\frac{r}{p}}}
		+ \pi \varepsilon
		\ll \varepsilon
	\end{flalign*}

	CASE (c$^\circ$): $\eta <m \leq \mu < M$. This can be split into two cases which are covered by Cases (a$^\circ$) and (b$^\circ$), respectively. Thus:
	\begin{equation*}
		\left\vert \sum\limits_{j=m}^{M} \sum\limits_{k=n}^{N} \right. c_{jk} \sin (jx) \sin(k\frac{2l_2\pi}{r})  \bigg\vert \ll \varepsilon .
	\end{equation*} We have two subcases to consider:

	Let $y = \frac{2l_2\pi}{r}, x \in  (\frac{(2l_1+1)\pi}{r},\frac{(2(l_1+1)\pi}{r})$ and let $\mu =  \left[\frac{1}{\frac{2(l_1+1)\pi}{r}-x}\right]^p$, where $l_1 \in \{  0, 1,...,[\frac{r}{2}]-1 \}$ and $ l_2 \in \{1,...,[\frac{r}{2}]-1 \}$, when $r\geq 4$ is an even number, and $l_1 \in \{  0, 1,...,[\frac{r}{2}]-1 \}$ and $ l_2 \in \{1,...,[\frac{r}{2}]\}$, when $r\geq 3$ is an odd number. Using  \eqref{eq_roznicowanie_druga_polowka} and \eqref{twierdzenie_lagrangea_drugi_przedzial} instead of \eqref{eq_roznicowanie_pierwsza_polowka} and \eqref{twierdzenie_lagrangea_pierwszy_przedzial} respectively, we can show similarly as in Case (a$^\circ$)-(b$^\circ$) that \eqref{zbieznosc_reszty_szeregu_podwojnego} holds.

	Similarly as above for $x=\frac{2l_1\pi}{r}$, $y \in (\frac{2l_2\pi}{r},\frac{(2l_2+1)\pi}{r}]$ or  $(\frac{(2l_2+1)\pi}{r},\frac{(2l_1+2)\pi}{r})$, using \eqref{dodatkowy_warunek_dla_max_k} instead of \eqref{dodatkowy_warunek_dla_max} we can show that \eqref{zbieznosc_reszty_szeregu_podwojnego} is satisfied.

	Summing up all the above estimations we, we conclude that \eqref{zbieznosc_reszty_szeregu_podwojnego} hold.
	This ends the proof. $\square$

	\subsection{Proof of Proposition \ref{proposition_rozbiezne_w_punkcie}}
	Let $p>1$, $c>1$ and for $m, n \in \mathbb{N}$ $c_{mn} = a_ma_n,\text{ where }$
	\begin{center}
		$a_n =\left\{ \begin{array}{l} \frac{3}{n\ln(n+1)},\text{ when } n=1\text{ (mod }3),
			\\\frac{1}{n\ln(n+1)},\text{  when } n=2\text{ (mod }3),
			\\\frac{1}{n\ln(n+1)},\text{  when } n=0\text{ (mod }3) \text{ and } n\neq 0\text{ (mod }6),
			\\\frac{1}{(n-3)\ln(n-2)} + \frac{1}{n^{1+\frac{1}{p}}\ln(n+1)},\text{  when } n=0\text{ (mod }6). \end{array} \right.$
	\end{center}
	First, we prove that $\{c_{mn}\} \in DGM(p,{}_2\alpha, {}_2\beta, {}_2\gamma, 3)$. Let
	\begin{align*}
		& A_n = \{ k \in \mathbb{N} : n\leq k \leq 2n-1\text{ and } k =1\text{ (mod }3) \}, \\
		&
		B_n = \{ k \in \mathbb{N} : n\leq k\leq 2n-1 \text{ and } k=1 \text{ (mod }3) \}, \\
		&
		C_n = \{ k \in \mathbb{N} : n\leq k\leq 2n-1\text{ and } k=0 \text{ (mod }3) \text{ and } k\neq 0\text{ (mod }6) \}. \\
		&
		D_n = \{ k \in \mathbb{N} : n\leq k\leq 2n-1\text{ and } k=0 \text{ (mod }6) \}.
	\end{align*}
	For $m,n \in \mathbb{N}$ we have:
	\begin{flalign*}
		\left( \sum_{k=n}^{2n-1} \right. & \left| c_{mk} - c_{m,k+3} \right|^p \Bigg)^{\frac{1}{p}}
		= \left( \sum_{k\in A_n} \left| c_{mk} - c_{m,k+3} \right|^p
		+ \sum_{k\in B_n} \left| c_{mk} - c_{m, k+3} \right|^p
		\right. & \\
		&
		\left.
		+ \sum_{k\in C_n} \left| c_{mk} - c_{m, k+3} \right|^p
		+ \sum_{k\in D_n} \left| c_{mk} - c_{m, k+3} \right|^p \right)^{\frac{1}{p}}
	\end{flalign*}
	It is clear that
	\begin{flalign*}
		& \left| \frac{1}{k\ln(k+1)} - \frac{1}{(k+3)\ln(k+4)} \right| \leq \frac{6}{k^2\ln(k+1)},
	\end{flalign*}
	Using elementary calculations and Lagrange's theorem, we have
	\begin{flalign*}
		&
		\left| \ln(k+4)-\ln(k-2) \right| \leq \frac{1}{k-2} \left( k+4-k+2 \right) = \frac{6}{k-2} &
	\end{flalign*}
	and
	\begin{flalign*}
		& \left| \frac{1}{(k-3)\ln(k-2)}  - \frac{1}{(k+3)\ln(k+4)} \right| 
		\leq
		\frac{(k-3)\left| \ln(k+4)-\ln(k-2)\right| + 6\ln(k+4)}{(k-3)(k+3)\ln(k-2)\ln(k+4)} &\\
		&
		\leq \frac{\frac{6(k-3)}{k-2}+ 6\ln(k+4)}{(k-3)(k+3)\ln(k-2)\ln(k+4)}
		\leq  \frac{12}{(k-3)(k+3)\ln(k-2)}
		\leq \frac{48}{k^2\ln(k+1)}.
	\end{flalign*}
	Hence, we have
	\begin{flalign*}
		\left( \sum_{k=n}^{2n-1} \right. & \left| c_{mk} - c_{m,k+3} \right|^p \Bigg)^{\frac{1}{p}}
		\leq \left( \sum_{k\in A_n} 3^p \left( \frac{6}{k^2\ln(k+1)} \right)^p \left\vert a_m \right\vert^p
		+ \sum_{k\in B_n} \left( \frac{6}{k^2\ln(k+1)} \right)^p \left\vert a_m \right\vert^p
		\right. \\
		&
		\left.
		+ \sum_{k\in C_n} \left( \frac{1}{k^{1+\frac{1}{p}}\ln(k+1)} \right)^p \left\vert a_m \right\vert^p
		+ \sum_{k\in D_n} \left( \frac{48}{k^2\ln(k+1)} + \frac{1}{k^{1+\frac{1}{p}}\ln(k+1)} \right)^p \left\vert a_m \right\vert^p
		\right)^\frac{1}{p} \\
		&
		\leq \frac{147}{n} \underset{b_2(n)\leq N\leq \lambda b_2(n)}{\text{max}} \sum\limits_{k=N}^{2N} \left\vert c_{mk} \right\vert ,
	\end{flalign*}
	where $ \lambda \geq 2$ and sequence $\left\{b_2(l)\right\}_{l=0}^{\infty}$ converges to infinity. In a similar way, we can show that
	\begin{flalign*}
		\left( \sum_{j=m}^{2m-1} \right. & \left| c_{jn} - c_{j+3,n} \right|^p \Bigg)^{\frac{1}{p}} \leq \frac{147}{n} \underset{b_1(m)\leq M\leq \lambda b_1(m)}{\text{max}} \sum\limits_{j=M}^{2M} \left\vert c_{jn} \right\vert, 
	\end{flalign*}
	where $ \lambda \geq 2$ and sequence $\left\{b_1(l)\right\}_{l=0}^{\infty}$ converges to infinity.\\
	Moreover we get
	\begin{flalign*}
		\left( \sum_{j=m}^{2m-1} \right. & \sum_{k=n}^{2n-1} \left\vert \Delta_{33} c_{jk}\right\vert^p \Bigg)^{\frac{1}{p}} =  & \\
		&
		\left( 
		\sum_{j\in A_{m}} \sum_{k \in A_{n}} \left\vert \Delta_{33} c_{jk}\right\vert^p
		+ \sum_{j\in B_{m}} \sum_{k \in A_{n}} \left\vert \Delta_{33} c_{jk}\right\vert^p
		+ \sum_{j\in C_{m}} \sum_{k \in A_{n}} \left\vert \Delta_{33} c_{jk}\right\vert^p
		+ \sum_{j\in D_{m}} \sum_{k \in A_{n}} \left\vert \Delta_{33} c_{jk}\right\vert^p 
		\right. \\
		&
		+ \sum_{j\in A_{m}} \sum_{k \in B_{n}} \left\vert \Delta_{33} c_{jk}\right\vert^p
		+ \sum_{j\in B_{m}} \sum_{k \in B_{n}} \left\vert \Delta_{33} c_{jk}\right\vert^p
		+ \sum_{j\in C_{m}} \sum_{k \in B_{n}} \left\vert \Delta_{33} c_{jk}\right\vert^p
		+ \sum_{j\in D_{m}} \sum_{k \in B_{n}} \left\vert \Delta_{33} c_{jk}\right\vert^p  \\
		&
		+ \sum_{j\in A_{m}} \sum_{k \in C_{n}} \left\vert \Delta_{33} c_{jk}\right\vert^p
		+ \sum_{j\in B_{m}} \sum_{k \in C_{n}} \left\vert \Delta_{33} c_{jk}\right\vert^p
		+ \sum_{j\in C_{m}} \sum_{k \in C_{n}} \left\vert \Delta_{33} c_{jk}\right\vert^p
		+ \sum_{j\in D_{m}} \sum_{k \in C_{n}} \left\vert \Delta_{33} c_{jk}\right\vert^p  \\
		&
		\left.
		+ \sum_{j\in A_{m}} \sum_{k \in D_{n}} \left\vert \Delta_{33} c_{jk}\right\vert^p
		+ \sum_{j\in B_{m}} \sum_{k \in D_{n}} \left\vert \Delta_{33} c_{jk}\right\vert^p
		+ \sum_{j\in C_{m}} \sum_{k \in D_{n}} \left\vert \Delta_{33} c_{jk}\right\vert^p
		+ \sum_{j\in D_{m}} \sum_{k \in D_{n}} \left\vert \Delta_{33} c_{jk}\right\vert^p \right)^{\frac{1}{p}} \\
		&
		\leq
		\frac{32 656}{mn} \sum_{j=m}^{2m-1} \sum_{k=n}^{2n-1} \left\vert a_ja_k\right\vert
		\leq {32 656}{mn} \underset{M+N\geq b_3(m+n)}{\text{sup}}\sum_{j=M}^{2M} \sum_{k=N}^{2N} \left\vert c_{jk}\right\vert.
	\end{flalign*}
	Hence $\left\{c_{mn}\right\} \in DGM(p,{}_2\alpha, {}_2\beta, {}_2\gamma,3)$.
	
	Now, we will show that $\left\{c_{mn}\right\} \notin DGM(1,{}_2\alpha, {}_2\beta, {}_2\gamma,3)$. For $m \in \mathbb{N}$ we have:
	\begin{flalign*}
		\sum_{k=n}^{2n-1} & \left\vert c_{mk} - c_{m,k+3} \right\vert
		\geq \sum_{k \in C_n} \left\vert c_{mk} - c_{m,k+3} \right\vert & \\
		&
		= \sum_{k \in C_n} \left\vert \frac{a_m}{k\ln(k+1)} - \left(\frac{1}{k\ln(k+1)} - \frac{1}{(k+3)^{1+\frac{1}{p}}\ln(k+4)} \right)a_m \right\vert \\
		&
		= \sum_{k \in C_n} \left\vert a_m \right\vert \frac{1}{(k+3)^{1+\frac{1}{p}}\ln(k+4)}
		\geq \frac{\left\vert a_m \right\vert}{(n+3)^{1+\frac{1}{p}}\ln(n+4)} \frac{n}{12}
		= \frac{\left\vert a_m \right\vert}{48(n+3)^{\frac{1}{p}}\ln(n+4)} 
	\end{flalign*}
	On the other hand, we get
	\begin{equation*}
		\frac{C}{n}  
		\underset{b_2(n)\leq N\leq \lambda b_2(n)}{\text{max}} \sum\limits_{k=N}^{2N} |c_{mk}| \leq C \frac{\left\vert a_m \right\vert}{n\ln(n+1)}.
	\end{equation*}
	Therefore, the inequality
	\begin{equation*}
		\sum_{k=n}^{2n-1} \left\vert c_{mk} - c_{m,k+3} \right\vert \leq {}_2\beta
	\end{equation*}
	can not be satisfied if $n \rightarrow \infty$.
	
	Now, we will show that the series \eqref{def_poswojny_szereg} is divergent in $(x_0, y_0) = (\frac{2}{3}\pi, \frac{2}{3}\pi)$. We have
	\begin{flalign*}
		\sum_{k=1}^{6N+5}& a_{k} \sin (kx_0) = 
		\sin(\frac{2}{3}\pi)\sum_{k=0}^{N} \left[ 
		\left( a_{6k+1} - a_{6k+2} \right) + \left( a_{6k+4} - a_{6k+5} \right)
		\right] & \\
		&
		= \sin(\frac{2}{3}\pi)\sum_{k=0}^{N} \left[ 
		\left( \frac{3}{(6k+1)\ln(6k+2)} - \frac{1}{(6k+2)\ln(6k+3)} \right)
		\right. \\
		&
		\left.
		+ \left( \frac{3}{(6k+4)\ln(6k+5)} - \frac{1}{(6k+5)\ln(6k+6)} \right)
		\right] \\
		&
		\geq 4\sin(\frac{2}{3} \pi) \sum_{k=0}^{N} \frac{1}{(6k+5)\ln(6k+5)} \rightarrow \infty \text{ as }N\rightarrow \infty.
	\end{flalign*}
	Thus the series \eqref{def_poswojny_szereg} is divergent in $(x_0, y_0) = (\frac{2}{3}\pi, \frac{2}{3}\pi)$. This ends our proof. $\square$	
\end{section}

\end{document}